\documentclass[11pt]{article}

\usepackage{amsmath}
\usepackage{amssymb}
\usepackage{epsfig}
\usepackage{graphicx}

\newtheorem{theorem}{Theorem}[section]
\newtheorem{lemma}{Lemma}[section]
\newtheorem{remark}[theorem]{Remark}

\newcommand{\qed}{\
\medskip\nopagebreak\nolinebreak\hfill$\square$\bigskip}

\newcommand{\be}{\begin{equation}}
\newcommand{\ee}{\end{equation}}
\newcommand{\beqa}{\begin{eqnarray}}
\newcommand{\eeqa}{\end{eqnarray}}

\newcommand{\NN}{\mathbb{N}}

\newcommand{\QQ}{\mathbb{Q}}
\newcommand{\RR}{\mathbb{R}}
\newcommand{\proof}{{\noindent \it Proof: }}

\newcommand{\Nabla}{\nabla\hspace{-2pt}}
\newcommand{\tfct}{\nu}

\newcommand{\R}{\mathbb{R}}

\def\e{\varepsilon}
\def\eps{\varepsilon}
\def\a{\'}
\def\g{\`}
\def\ds{\displaystyle}
\def\div{\nabla\cdot}

\def\kkvec{{\mathbf K}}
\def\xvec{{\mathbf x}}
\def\uvec{{\mathbf u}}
\def\lvec{{\mathbf L}}
\def\jvec{{\mathbf J}}
\def\yvec{{\mathbf y}}
\def\vvec{{\mathbf v}}
\def\evec{{\mathbf E}}
\def\evecr {{{\mathbf E}_r}}
\def\gvec{{\mathbf G}}
\def\avec{{\mathbf A}^\varepsilon}
\def\eext{{\mathbf \Xi}}
\def\eextr{{{\mathbf \Xi}_r}}

\def\cale{{\mathcal E}}
\def\caler{{\mathcal E}_r}
\def\calo{{\mathcal O}}
\def\calq{{\mathcal Q}}

\begin{document}

\title{Long time simulation of a beam in a 
periodic focusing channel via a two-scale PIC-method}

\bibliographystyle{plain}

\author{
E. Fr\'enod\thanks{LMAM et Lemel, Universit\'e de Bretagne sud,
Campus de Tohannic, 56000 Vannes, France} 
\and F. Salvarani\thanks{
Dipartimento di Matematica,
Universit\g a degli Studi di Pavia.
Via Ferrata, 1 - 27100 Pavia, Italy}
\and E. Sonnendr\"ucker\thanks{Institut de Recherche Math\'ematique
Avanc\'ee UMR 7501, Universit\'e Louis Pasteur, Strasbourg 1 et CNRS, et Projet Calvi INRIA Nancy - Grand Est,
7 rue René Descartes - 67084
Strasbourg Cedex, France}
}

\date{\today}
\maketitle

\begin{abstract}
We study the two-scale asymptotics for a charged beam under
the action of a rapidly oscillating external electric field.
After proving the convergence to the correct asymptotic state, we
develop a numerical method for solving the limit model involving two time scales
and validate its efficiency for the 
simulation of long time beam evolution.
\end{abstract}



\section{Introduction}
\label{intro}

Charged particle beams, as used in beam physics for many applications ranging from heavy ion fusion 
to cancer therapy, often need to be transported on very long distances in particle accelerators. 
Moreover, as particles of the same charge repulse each other, these beams need to be focused using external 
electric or magnetic fields (see, for example, the paper of Filbet and Sonnendr\"ucker 
\cite{FilSon} for a mathematical description of particle beam modelling and numerical simulation).
Beam focusing can be performed either in linear or circular accelerators using a periodic focusing lattice.

The aim of this paper is to study the evolution of a beam over a large number of periods.
The motion of electrically charged particles is governed by several phenomena: {\it in primis}, it is
necessary to take into account the interactions between the electric fields generated by the
particles themselves and by the external focusing electromagnetic field.
The kinetic approach is a successful method suitable to model a system composed by a great
number of particles; in our case the modelling framework will be
given by coupling a kinetic equation with the Maxwell equation.
When it is possible to disregard the collisions between particles, the kinetic part of the
modelling analysis is performed by means of the Vlasov collisionless equation \cite{Gla}.
In this paper we will consider only non-relativistic long and thin beams; therefore,
instead of studying the phenomenon by means of the full Vlasov-Maxwell system,
we will consider its paraxial approximation.
This simplified model of the Vlasov-Maxwell 
equations is particularly adapted to the study of long and thin beams and describes
the evolution of a charged particle beam in the plane transverse 
to its direction of propagation.

The key assumptions of the paraxial approximation are the following:
\begin{enumerate}
\item we suppose that there is a predominant length scale, namely the longitudinal length
of the beam, whereas the transverse thickness of the beam is negligible with
respect to its longitudinal length;
\item we will assume that the beam has already reached its stationary state.
\end{enumerate}
Such assumptions are often verified in applications: in a standard accelerator,
the orders of magnitude of the different lengths taken into account verify the
first assumption;
moreover it is often interesting to study stable (in time) configurations of the 
beam, and therefore the second assumption applies.
In the paraxial framework, the effects of the self consistent magnetic and electric
fields can be both taken into account by solving a single Poisson equation.
Hence the model becomes similar to the Vlasov-Poisson system in the transverse plane
with respect to the accelerator axis. In the paraxial model, the $z$ coordinate plays the role
of time. We will switch here to more conventional notations and consider the time-dependent
two-dimensional Vlasov-Poisson system.
We assume moreover, in accordance with most experiments, that 
the average external field is large compared to the
self-consistent field.
More details about the derivation of this model can be found in \cite{DegRav} or
\cite{FilSon}, whereas a more complete description on accelerators' physics
can be found in the book by Davidson and Qin \cite{DavQin}.

In the framework we just precised, the Vlasov-Poisson system
takes the form
\begin{equation}
\label{VP-paraxial}
\left \{
\begin{array}{llll}
\displaystyle
\frac{\partial f_\e}{\partial t} + \frac 1 \e \vvec\cdot\Nabla_x f_\e
+\left( \evec_\e+ \eext^\e\right ) \cdot\Nabla_v f_\e=0,
\\
\\
\displaystyle
\evec_\e= -\nabla \Phi_\e,
\\
\\
\displaystyle
\Delta \Phi_\e= -\rho_\e(t, \xvec), \qquad \rho_\e(t,
\xvec)=\int_{\R^2_v}f_\e(t,\xvec,\vvec)\, d\vvec,
\\
\\
\displaystyle
f_\e(t=0,\xvec,\vvec)= f_0,
\end{array}
\right.
\end{equation}
where $f_\e=
f_\e ( t,\xvec,\vvec)$ is the distribution function, $t\in [0,T)$ for some $T<\infty$,
$\xvec=(x_1,x_2)\in\R^2_x$ is the position vector and
$\vvec=(v_1,v_2)\in\R^2_v$ is the velocity vector.
We will denote by $\Omega=\R^2_x\times\R^2_v$.

As said before, we are working in a stationary setting.
In this context, the variable $t$ does not represent, from a physical point of view, a
time variable, but rather the longitudinal coordinate. We have nevertheless chosen
to use the variable $t$ because of the great similarity of the paraxial approximation
with respect to a two-dimensional (in space) time-dependent Vlasov equation.
The electric field $\eext^\e(t,\xvec)$ is supposed to be external, and has
the form
$$
\eext^\e(t,\xvec)= 
-\frac{1}{\e}H_0 \xvec + H_1 \hspace{-2pt}\left(\omega_1 \frac{t}\e\right) \kkvec \xvec,
$$
where $H_0$ is a positive constant and $H_1$ is $2\pi-$periodic,
with mean value $0$.
$\kkvec$ is a $2\times 2$ symmetrical isometric matrix. 
In the two most common physical situations, $\kkvec=I$ or 
$\kkvec=\begin{pmatrix} 1&0\\0&-1 \end{pmatrix}$.
The parameter $\e$ is a scaling parameter acting on the "time" scale 
(i.e. the longitudinal scale). Its physical meaning will 
be detailed in the next section.

The form of the external electric field will satisfy some relevant physical
requirements.
The field which models the periodic focusing lattice needs indeed to be
linear in $\xvec$ to compensate for the self field which is also linear in $\xvec$
for a uniform beam but with the opposite sign. The variation of the field
with respect to $t$ is assumed to be a 
small amplitude field oscillating at the same period as the lattice.

\smallskip
Under additional physical hypotheses, a simplified version of (\ref{VP-paraxial})
may be considered. These assumptions consist, in the case when $\kkvec=I$, in considering
that the beam is axisymmetric, \emph{i.e.} invariant under rotation in $\R^2_x$.
Then, by writing (\ref{VP-paraxial}.a) in polar coordinates $(r,\theta)$, such that
$r=|\xvec|$, $x_1= r \cos(\theta)$ and $x_2= r \sin(\theta)$ we get 
\begin{equation}
\frac{\partial f_\e}{\partial t} 
+ \frac 1\e v_r \frac{\partial f_\e}{\partial r} 
+\left(\evecr_\e+\eextr^\e+ \frac{(r v_\theta)^2}{r^3} \right)\frac{\partial f_\e}{\partial v_r} 
=0,
\end{equation}
where $v_r$ is the projection of the velocity on radius $\xvec$. In other words, 
$v_r=\vvec\cdot \xvec/r$ and $v_\theta=\vvec\cdot \xvec^\perp/r$ where 
$\xvec^\perp=(-x_2,x_1)$.
If finally, $f_\e$ is assumed to be concentrated in angular momentum $r v_\theta=0$
(this property is achieved if the initial condition $f_0$ is concentrated in 
angular momentum), system (\ref{VP-paraxial}) yields
\begin{equation}
\label{2}
\left \{
\begin{array}{llll}
\displaystyle
\frac{\partial f_\e}{\partial t} 
+ \frac 1\e v_r \frac{\partial f_\e}{\partial r} 
+\left(\evecr_\e+\eextr^\e\right)\frac{\partial f_\e}{\partial v_r} 
=0,
\\
\\
\displaystyle
\frac 1r \frac{\partial (r\evecr_\e)}{\partial r} 
= \rho_\e(t, r), \qquad 
\rho_\e(t, r)=\int_{\R}f_\e(t,r,v_r)\, d v_r,
\\
\\
\displaystyle
f_\e(t=0,r,v_r)= f_0,
\end{array}
\right.
\end{equation}
where $f_\e=f_\e(t,r,v_r)$ for $t\in[0,T)$, $r\in \R^+$
and $v_r\in \R$. The external force $\eextr^\e$ writes here
$$
\eextr^\e(t,r)= 
\big(-\frac{1}{\e}H_0 + H_1 \hspace{-2pt}\left(\omega_1 \frac{t}\e\right)\big) r.
$$

System (\ref{2}) is naturally set for $r\in \R^+$. Nevertheless, we can consider 
that $r\in\R$ using the convention  $f_\e(t,r,v_r)=f_\e(t,-r,-v_r)$ and
$\eextr^\e(t,r)=\eextr^\e(t,-r)$. This is what we do in the following.

\bigskip
The main purpose of the paper is the study of systems (\ref{VP-paraxial})
and (\ref{2}) in the limit as $\e\to 0$.
We have chosen to work in the framework of the so-called two-scale convergence, an
homogenization technique introduced by N'guetseng
\cite{Ngu} and subsequently developed by Allaire \cite{All}
and  Fr\'enod, Raviart and Sonnendr\"ucker \cite{FreRavSon},
which essentially says that $f_\e(t,\xvec,\vvec)$ is close to
$F(t,t/\e,\xvec,\vvec)$ for a given profile $F$.
This notion is stronger than the usual weak or weak* convergence and
is very adapted
to the study of asymptotic behaviour of transport equations with rapidly 
oscillating terms, as was shown in \cite{FreTh, FreHam, FreRavSon, FreSon0,FreSon}.
On  the other hand, we develop a numerical method in order to compute an approximation of $f_\e$. 
This method consists, in the spirit suggested in  Fr\' enod, Raviart and Sonnendr\" ucker 
\cite{FreRavSon} and Fr{\'e}nod \cite{FreSpds} and explored in Ailliot, Fr\'enod and Monbet \cite{AilFreMon} and in Fr\'enod, Mouton and Sonnendr\"ucker \cite{FreMouSon},
in computing the profile $F$ by discretizing the equation it satisfies.
Then an approximation of $f_\e$ is reconstructed using 
$f_\e(t,\xvec,\vvec) \sim F(t,t/\e,\xvec,\vvec)$.
For simplicity reasons, the considered method is implemented in the simplified 
case of system (\ref{2}). Nevertheless, it is easy to see that it may be adapted in a
straightforward way to system (\ref{VP-paraxial}) and others situations where multiple time
scales are involved.

The paper is organized as follows: in Section \ref{scaling} we obtain the scaled
Vlasov-Poisson system (\ref{VP-paraxial}) from the standard Vlasov-Poisson system 
by means of a scaling procedure.
Section \ref{well} is then devoted to the derivation of some useful properties
of the solution of (\ref{VP-paraxial}), whereas Section \ref{homogeneization} is
dedicated to homogenization of systems (\ref{VP-paraxial}) and (\ref{2}) by means of 
the two-scale convergence.
Finally, in Section \ref{numerics} we will implement and test the previously 
described numerical method based on the homogenization analysis.

\section{Scaling of the Vlasov equation: the paraxial approximation}
\label{scaling}

We present in this paragraph the scaling procedure leading to
Equation (\ref{VP-paraxial}).
We first note that, in the formulation of the problem, there are different
scales: the length of the accelerator $L$, the transverse radius of the
beam $\lambda$ and the period of the oscillating external field $l$.
In a standard accelerator, the orders of magnitude of such quantities are:
$L\simeq 10^3$ m; $l\simeq 10^{-1}$ m and, finally, $\lambda\simeq 10^{-2}$
m.
The standard three-dimensional Vlasov-Poisson system, in presence of an
external electric field
$\eext$, is given by
\begin{equation}
\label{VP}
\left \{
\begin{array}{llll}
\displaystyle
\frac{\partial f}{\partial t} +\vvec\cdot\Nabla_x f
+\frac{q}{m} \left( \evec+ \eext\right ) \cdot\Nabla_v
f=0,
\\
\\
\displaystyle
\evec= -\nabla \Phi,
\\
\\
\displaystyle
\Delta \Phi= -\frac 1{\e_0}\rho(t,\xvec), \qquad \rho(
\xvec,t)=q\int_{\R^3_v}f(t,\xvec,\vvec)\, d\vvec,
\\
\\
\displaystyle
f(t=0,\xvec,\vvec)= f_0,
\end{array}
\right.
\end{equation}
where $f\equiv
f(t,\xvec,\vvec)$ is the distribution function, $t\in [0,T)$ for some $T<\infty$,
$\xvec=(x_1,x_2,x_3)\in\R^3_x$ is the position vector and
$\vvec=(v_1,v_2,v_3)\in\R^3_v$ is the velocity vector.

The external electric field $\eext=
\eext(\xvec)$ is supposed to be independent of time and periodic with
respect to the variable $x_3$.
Under the paraxial approximation, $f$ is supposed to be stationary
with respect to time and the velocity of the particles
with respect to the longitudinal direction of the beam (which will
be ${\mathbf e}_3$ throughout the whole paper) is constant, that
is $v_3=v_b$.
These assumptions will enable us to eliminate the time derivative in the
first equation of (\ref{VP}) and suppose that $f(t,\xvec,\vvec)=
f(\xvec,v_1, v_2) \otimes \delta(v_b-v_3)$.
In order to obtain the dimensionless version of the Vlasov-Poisson
system
(\ref{VP}), we define the new variables $x_i'$, $v_i'$ by:
$$
x_i=\lambda\, x_i', \qquad v_i=v_b\, v_i', \qquad i=1,\, 2
$$
$$
x_3=L\, x_3', \qquad v_3=v_b.
$$
Moreover, the dimensionless density $f'$, external electric field $\eext'$,
self-consistent electric field  $\evec'$
and its potential $\Phi'$ are defined by
$$
f=\bar f f',
\qquad  \eext= \bar \Xi\, \eext', 
\qquad  \evec=\bar E \, \evec', 
\qquad \Phi=\bar\Phi\, \Phi'.
$$
Under such a change of variables, we obtain that the Vlasov-Poisson system (\ref{VP})
can be written in the following form:
\begin{equation}
\left \{
\begin{array}{llll}
\displaystyle
\frac {\bar f v_b}L
\frac{\partial f'}{\partial x_3'} + \frac {\bar f v_b}\lambda \sum_{i=1}^2 v_i'
\frac{\partial f'}{\partial x_i'}
+ \frac {q \bar E \bar f }{mv_b}
\sum_{i=1}^2 \left( E_i'+\frac{\bar \Xi}{\bar E} \Xi_i'\right
)\frac{\partial f'}{\partial v_i'} =0,
\\
\\
\displaystyle
\bar E E_i'= - \frac {\bar \Phi}\lambda \frac{\partial \Phi'}{\partial x_i'}, 
\qquad i=1, \ 2; 
\qquad \bar E E_3'= - \frac {\bar \Phi}L \frac{\partial \Phi' }{\partial x_3'},
\\
\\
\displaystyle
\frac{\bar \Phi}{\lambda^2}
\sum_{i=1}^2 \left  (
\frac{\partial^2 \Phi'}{\partial x_i^{\prime 2}} \right )
+ \frac{\bar \Phi}{L^2} \frac{\partial^2 \Phi'}{\partial x_3^{\prime 2}}
= - \rho'(\xvec'), \qquad \rho'(
\xvec')=\frac{v_b^2  \bar f q}{\e_0} \int_{\R^2_{v'}}f'(\xvec',\vvec')\, d\vvec'
\\
\\
\displaystyle
\bar f f'(t=0,\xvec',\vvec')= f_0(\lambda x_1',\lambda x_2', L x_3', v_b v_1, v_b v_2 ).
\end{array}
\right.
\end{equation}
It is natural to choose, as unit measure of the electric field
$$
\bar E=\frac{mv_b^2}{qL};
$$
moreover, since the third component of the electric field will
play no role, its unity measure induces that
$$
\qquad \bar\Phi=\bar E \lambda= \frac{m\lambda v_b^2}{qL}.
$$
Finally, if we assume to have a small self consistent electric field,
we are forced to work with a small density of particles. Therefore, we
will assume that the dimensionless form of the density $f$ is given by
$f=\bar f\, f'$, with
$$
\bar f=\frac{m \e_0}{q^2\lambda L}.
$$

As mentioned in the introduction, we consider here low intensity beams
so that the external electric field is much stronger than
the self-consistent field. 
We translate this, introducing a small parameter $\e$, by saying that
$\bar \Xi / \bar E = 1/\e$ and assuming that $\eext'$ writes
$$
\eext'= 
\left [H_0\gvec_1 + 
\e  H_1\hspace{-2pt}\left ( \omega_1\frac{Lx_3'}l\right )\gvec_2 \right ],
$$
where $H_1$ is $2\pi-$periodic and for a $2\times 2$ isometric matrix $\kkvec$
$$
\gvec_1= - \begin{pmatrix}x_1\\x_2\\ 0\end{pmatrix}
\text{ ~ and ~ }
\gvec_2= \begin{pmatrix} \kkvec&0\\0&0\end{pmatrix}
\begin{pmatrix}x_1\\x_2\\ 0\end{pmatrix}.
$$
On the other hand, we also suppose that the ratio between the average transverse
radius of the beam and the longitudinal length of the accelerator $\lambda/L=\e$,
and, that the ratio between the period of the external electric field and the
macroscopic quantity $l/L=\e$.
Forgetting the terms in $\e^2$, 
we have then obtained the dimensionless Vlasov-Poisson system in the
paraxial approximation:
\begin{equation}
\label{VP-paraxial-primed}
\left \{
\begin{array}{llll}
\displaystyle
\frac{\partial f'}{\partial x_3'} + \frac 1 \e \sum_{i=1}^2 v_i'
\frac{\partial f'}{\partial x_i'}
+\sum_{i=1}^2 \left( E'_i+ \frac1\e\eext'\right)
\frac{\partial f'}{\partial v_i'} =0,
\\
\\
\displaystyle
E'_i= -\partial \Phi' /\partial x_i', \qquad i=1,2,
\\
\\
\displaystyle
\sum_{i=1}^2 \frac{\partial^2 \Phi'}{\partial x_i^{\prime 2}}= - \rho'(
\xvec'), \qquad \rho'(
\xvec')=\int_{\R^2_{v'}}f'(\xvec',\vvec')\, d\vvec',
\\
\\
\displaystyle
f'(t=0,\xvec',\vvec')= f_0',
\end{array}
\right.
\end{equation}
with
$$ 
\frac{1}{\e} \eext'= -
\frac{1}{\e} H_0 \gvec_1 + 
H_1\hspace{-2pt} \left (\omega_1 \frac{x_3'}\e \right ) \gvec_2 .
$$
Denoting the variable $x_3'$ with $t$ (in order to have in (\ref{VP-paraxial-primed})
only two-dimensional vectors)
and eliminating all the primed variables,
we finally obtain the scaled Vlasov-Poisson paraxial system
(\ref{VP-paraxial}).

\section{Main properties of the solution}
\label{well}

We recall that, from now on, all the vectors are intended to be two-dimensional.

In the sequel, we will use heavily the property of boundedness of
the solution with respect to some norms. This is guaranteed by the
following lemma:

\begin{lemma}
\label{bounds}
Let $f_\eps$ be the solution of the Vlasov-Poisson system $(\ref{VP-paraxial})$,
with non-negative (a.e.) initial data $f_0$ of class $(L^1\cap L^p)(\Omega)$, 
$p\geq 2$ and such that the moment of order two is finite:
$$
\int_\Omega (|\xvec|^2 + |\vvec|^2 )f_0\, d\xvec \, d\vvec
< + \infty.
$$
Then $f_\eps$ is bounded in $L^\infty([0,T];L^p(\Omega))$ uniformly with respect to 
$t$. Moreover 
\begin{equation}\label{e1}
\Vert (|\xvec|^2 f_\e+ |\vvec|^2 f_\e) \Vert_{L^\infty([0,T];L^1
(\Omega))} \leq C \, e^{\alpha T},
\end{equation}
with
$$
\alpha=\frac {\eps \Vert H_1\Vert_\infty}{\min \{1\, , \, H_0 \} }.
$$
Finally
\begin{equation}\label{e2}
\|\rho_\eps(\xvec,t)\|_{L^\infty([0,T];L^{\frac{3}{2}}(\R_x^3))}\leq C e^{\alpha T/3},
\end{equation}
for some constant $C$.
\end{lemma}

\proof
We multiply the Vlasov equation (\ref{VP-paraxial}a) by $f_\eps^{p-1}$ and integrate
in $\xvec$ and $\vvec$. We obtain
$$
\| f^\eps\|_{L^\infty([0,T];L^p(\Omega))} \leq C,
$$
for some constant $C$.
This proves the first part of the lemma.

In order to prove the other statements, we
multiply the Vlasov equation (\ref{VP-paraxial}a)
by $|\vvec|^2$, and integrate with respect to
$\xvec$ and $\vvec$.
We get
\begin{equation}
\label{vlasen}
\frac{d}{dt}\int_\Omega f_\eps |\vvec|^2 \, d\xvec \,d\vvec
-2\int_{\R^2_x}\jvec\cdot(\evec_\e + \eext^\e)\, d\xvec =0,
\end{equation}
where 
$$
\jvec(\xvec,t)=\int_{\R^2_v} \vvec f_\eps\,d\vvec.
$$
Now, integrating the Vlasov equation in $\vvec$ gives
the continuity equation
\begin{equation}
\label{conteq}
\frac{\partial \rho_\e}{\partial t}+\frac 1\eps \div\jvec=0.
\end{equation}
Thanks to the continuity equation (\ref{conteq}), we note that
\begin{equation}
\label{truc1}
\int_{\R^2_x} \jvec\cdot\evec_\e\,d\xvec = 
-\int_{\R^2_x} \jvec\cdot\nabla \Phi_\e\,d\xvec = 
\int_{\R^2_x} \div\jvec\, \Phi_\e\,d\xvec
=  -\eps\int_{\R^2_x}\frac{\partial \rho_\e}{\partial t}\,\Phi_\e\, d\xvec.
\end{equation}
Using now the Poisson equation, we get
\begin{equation}
\label{truc2}
\frac{1}{2}\frac{d}{dt}\int_{\R^2_x}(\nabla \Phi_\e)^2\,d\xvec = 
-\int_{\R^2_x}\frac{\partial }{\partial t} \left( \Delta\Phi_\e\right) \,\Phi_\e\,d\xvec =
\int_{\R^2_x} \frac{\partial\rho_\e}{\partial t}\,\Phi_\e\,d\xvec.
\end{equation}
Coupling (\ref{truc1}) and (\ref{truc2}) we deduce:
\begin{equation}
\label{truc12}
-2\int_{\R^2_x}\jvec\cdot\evec_\e\,d\xvec=\eps 
\frac{d}{dt}\int_{\R^2_x}(\nabla \Phi_\e)^2\,d\xvec.
\end{equation}

On the other hand, the external electric field $\eext^\e$ can be
also derived from the potential
$$
P(\xvec)= \frac{1}{2\e}H_0|\xvec|^2
- \frac{1}{2}H_1\hspace{-2pt}\left ( \omega_1\frac{t}\e\right )\kkvec\xvec\cdot\xvec . 
$$
A variant of the procedure used for the self-consistent field
allows us to deduce
$$
-2 \int_{\R^2_x} \jvec\cdot\eext^\e\,d\xvec = 
2 \int_{\R^2_x} \jvec\cdot\nabla P\,d\xvec = 
-2\int_{\R^2_x} \div\jvec\, P\,d\xvec
=  2\eps\int_{\R^2_x}\frac{\partial \rho_\e}{\partial t}\,P\, d\xvec =
$$
$$
H_0 \int_{\R^2_x}|\xvec|^2 \frac{\partial \rho_\e}{\partial t}\, d\xvec-
\eps\int_{\R^2_x}  H_1\hspace{-2pt}\left (\omega_1 \frac{t}\e\right )\kkvec\xvec\cdot\xvec 
\frac{\partial \rho_\e}{\partial t}\,d\xvec,
$$
and therefore (\ref{vlasen}) becomes
$$
\frac{d}{dt}\biggl[
\int f_\eps |\vvec|^2 \,d\vvec\,d\xvec + \eps\int(\nabla \Phi_\e)^2\,d\xvec
+ H_0\int_{\R^2_x}|\xvec|^2  \rho_\e\, d\xvec
\biggr]
=
$$
$$
=
\eps H_1\hspace{-2pt}\left (\omega_1 \frac{t}\e\right )\int_{\R^2_x} \kkvec\xvec\cdot\xvec 
\frac{\partial \rho}{\partial t}\,d\xvec.
$$
If we multiply the Vlasov Equation (\ref{VP-paraxial}a) by $|\xvec|^2$ and then
integrate with respect to $\xvec$ and $\vvec$, we deduce that
$$
\int_{\R^2_x} \kkvec\xvec\cdot\xvec  \,\frac{\partial \rho_\e}{\partial t}\,d\xvec
=- \int_{\Omega} \kkvec\xvec\cdot\xvec \, \vvec \cdot\Nabla_x f_\e \,d\xvec\, d\vvec
= 2 \int_{\Omega} (\kkvec\xvec \cdot \vvec)\, f_\e \,d\xvec\, d\vvec
$$
and thanks to the elementary inequality 
$2\kkvec\xvec\cdot\vvec \leq  (|\kkvec\xvec|^2+|\vvec|^2)\leq  (|\xvec|^2+|\vvec|^2)$,
we obtain that
$$
\frac{d}{dt}\biggl[
\int_{\Omega} f_\eps |\vvec|^2 \,d\vvec\,d\xvec + \eps\int_{\R^2_x}(\nabla \Phi_\e)^2\,d\xvec
+ H_0 \int_{\R^2_x}|\xvec|^2 \rho_\e \, d\xvec
\biggr]
\leq 
$$
$$
\leq
\eps \Vert H_1\Vert_\infty \biggl[
\int_{\Omega} |\xvec|^2 f_\e \,d\xvec\, d\vvec + 
\int_{\Omega} |\vvec|^2 f_\e \,d\xvec\, d\vvec \biggr].
$$
Finally
$$
\Vert (|\xvec|^2 f_\e+ |\vvec|^2 f_\e) \Vert_{L^\infty([0,T];L^1
(\Omega))}
\leq C \, e^{\alpha T},
$$
with
$$
\alpha=\frac {\eps \Vert H_1\Vert_\infty}{\min \{1\, , \, H_0 \} },
$$
which proves the second statement of the lemma.

The proof of the bound on $\rho_\e$ is a straightforward
consequence of the following classical estimate 
(see for instance \cite{FreSon} (lemma 4.4)):
$$
\int_{\R^2_x} |\rho_\e(\xvec,t)|^{3/2}\,d\xvec \leq C_3 
\left (\int_{\Omega} (f_\eps)^2\,d\xvec\,d\vvec\right )^{1/2} \left(\int_\Omega
|v|^2 f_\eps\,d\xvec\,d\vvec\right)^{1/2}.
$$
Hence the proof of the lemma is  complete.
\qed

\section{Homogenization of the paraxial approximation}
\label{homogeneization}
Our goal consists in deducing the equations satisfied by the limit of
$f_\e$ as $\e\to 0$.
The homogenization of partial differential equations with rapidly
periodic coefficients, as in the present case, can be fruitfully studied 
applying the two-scale convergence, introduced by N'guetseng
\cite{Ngu} and subsequently developed by Allaire \cite{All}
and  Fr\'enod, Raviart and Sonnendr\"ucker \cite{FreRavSon}
and which essentially stands in the following theorem.

\begin{theorem}
\label{Allaire}
Let $\calq$ be a regular subset of $\R^n$, $X$ a Banach space and $X'$
its dual space, with $\left <\;,\;\right >$ as duality bracket. Let $1<p\leq+\infty$
and $p'$ be such that $\frac1p+\frac1{p'}=1$ and let 
$Y=[0,a_1]\times\dots\times[0,a_n]$ defined for finite real numbers $a_i$.
\newline
If a sequence $(u_\e) = (u_\e(q))$ is bounded in $L^{p'}(\calq;X')$
then there exists a function $u_0=u_0(q,y) \in L^{p'}(\calq\times Y;X')$
such that,  up to a subsequence,
$$
\lim_{\e\to 0}  \int_\calq \big< u_\e(q),\tfct\left (q,\frac q\e \right )\big>\, dq 
=
\int_\calq \int_Y \big< u_0(q,y),\tfct\left (q,y\right )\big>\, dq dy,
$$
for any function $\tfct \in L^{p}(\calq\, ; \, C_Y(\R^n;X))$,
where $ C_Y(\R^n;X)$ stands for the space of continuous
$Y-$periodic functions on $\R^n$ with values in $X$.
We say that $(u_\e(x))$ two-scale converges to $u_0(x,y)$.
\end{theorem}

More precisely, we will use a theorem of Fr\a enod and Sonnendr\" ucker
\cite{FreSon}, which gives a generic framework, in
the context of two-scale convergence, for linearly
perturbed conservation laws of the form
\begin{equation}
\label{abseq2}
\left \{
\begin{array}{ll}
\ds
\frac{\partial u_\e}{\partial t}+\avec\cdot\Nabla_x u_\e +
\frac{1}{\e}\lvec\cdot\Nabla_x u_\e = 0, 
\\
\\ 
\ds u_\e(t=0)=u_0.
\end{array}
\right.
\end{equation}
In this system, $u^{\eps}\equiv u^{\eps}( t,\xvec),$ $t\in[0,T)$ for
some $T<\infty$ and $\xvec\in\mathbb{R}^n=\calo$. Moreover
it is assumed that, for all $\e>0$,
$\Nabla_x\cdot\avec=0$, and that, for some $q>1$, 
$\avec\equiv\avec(t,\xvec)$ two-scale converges to
${\mathcal A}\equiv{\mathcal A} (t,\tau,\xvec)\in{L^\infty}
([0,T]\times[0,\theta]\,; \, W^{1,q}(K))$
for all compact sets $K\in \R^n$.
Finally, ${\mathbf L} = M\xvec $ where $M$ is
a real $n\times n$ matrix with constant entries, satisfying tr$M =0$ and
such that $e^{\tau M}$ is $\theta-$periodic.

Under these conditions,  the following theorem is proved in \cite{FreSon}:

\begin{theorem}
\label{thcali}
Under the assumptions above, if moreover the sequence $(u_\eps)$ of
solution of $(\ref{abseq2})$ satisfies
\begin{equation}\label{uepsborne1}
\|u_\eps\|_{L^\infty([0,T];L^p(\calo))} \leq C,
\end{equation}
for some $p>1$ such that $1/p+1/q'<1$, where $1/q'= \max \{   1/q-1/n\, ,\, 0  \}$,
then, extracting a subsequence, 
$u_\eps$ two-scale converges to a profile 
$$
U\in L^\infty([0,T]\times[0,\theta];L^p(\calo)).
$$
Moreover we have
\begin{equation}
\label{UtoU0}
U(t,\tau,\xvec )= U_0(t,e^{-\tau M}\xvec),
\end{equation}
where $U_0\equiv U_0(t,\yvec)$ 
is solution to
\begin{equation}
\label{eqU0}
\left \{
\begin{array}{ll}\ds
\frac{\partial U_0}{\partial t} + \int_0^\theta e^{-\sigma M}
{\mathcal A} (t,\sigma,e^{\sigma M}\yvec)\, d\sigma\cdot\Nabla_y U_0=0,
\\
\\
\ds
{U_0}(t=0) =\frac{1}{\theta} u_0.
\end{array}
\right.
\end{equation}
\end{theorem}
The problem of homogenizing system (\ref{VP-paraxial}) enters the framework we just
presented with $\xvec$ replaced with $(\xvec,\vvec)\in \Omega$,
with 
\begin{equation}\label{Aeps} 
\avec = 
\begin{pmatrix}
0 \\ \evec _\e(t, \xvec)+H_1\hspace{-2pt}\left (\omega_1 \frac{t}\e\right )\kkvec \xvec
\end{pmatrix},
\end{equation}
and with
$$\lvec=\begin{pmatrix} \vvec \\ -H_0 \xvec\end{pmatrix}
\text{ or in other words } M =\begin{pmatrix} 0&I\\ -H_0 I &0 \end{pmatrix}.$$
It is an easy game to see that
\begin{equation}\label{exptauM}
e^{\tau M}=
\left( \begin{array}{cc}
{\mathcal R}_1(\tau) & {\mathcal R}_2(\tau) /\omega_0 
\\
-\omega_0 {\mathcal R}_2(\tau) & {\mathcal R}_1(\tau)
\end{array} \right),
\text{ with }\omega_0 =\sqrt{H_0},
\end{equation}
where
$$
{\mathcal R}_1(s)=
\left( 
\begin{array}{cc}
\cos\left (\omega_0 s \right ) & 0
\\
0 & \cos\left (\omega_0 s \right ) 
\end{array} \right),
$$
and
$$
{\mathcal R}_2(s)=
\left( \begin{array}{cc}
\sin\left (\omega_0 s \right )& 0 
\\
0 & \sin\left (\omega_0s \right )
\end{array} \right).
$$
We can now prove:
\begin{theorem}
\label{main}
Let $f_0$ satisfy the hypotheses of Theorem $\ref{bounds}$. Then, if
we consider a sequence of solutions $(f_\e,\evec_\e)$ depending on $\e$,
extracting a subsequence, we have that $f_\e$ two-scale converges to $F
\in L^\infty([0,T]\times[0,\frac{2\pi}{\omega_0}];L^2(\Omega))$ 
and $\evec_\e$ two-scale converges to 
$\cale \in L^\infty([0,T]\times[0,\frac{2\pi}{\omega_0}];W^{1,3/2}(\R^2_x))$.
\newline
Moreover, there exists a function 
$G=G(t, \yvec, \uvec)\in L^\infty([0,T];L^2(\Omega))$ 
such that
$$
F(t,\tau, \xvec,\vvec)= G \left( t,{\mathcal R}_1(-\tau)\xvec +
\frac1{\omega_0}{\mathcal R}_2(-\tau)\vvec 
\, , \, - \omega_0{\mathcal R}_2(-\tau)\xvec+
{\mathcal R}_1(-\tau)\vvec \right ).
$$
The pair $(G,\cale)$ is finally the solution of 
\begin{equation}
\label{eqG}
\left \{
\begin{array}{ll}\ds
\frac{\partial G}{\partial t} + \frac{1}{\omega_0}
\int_0^{2\pi/\omega_0} {\mathcal R}_2(-\sigma)
\cale (t,\sigma,{\mathcal R}_1(\sigma)\yvec+
\frac{1}{\omega_0}{\mathcal R}_2(\sigma)\uvec)\, d\sigma
\cdot\Nabla_y \, G
\\ \ds~~~~~~~~~
+\int_0^{2\pi/\omega_0} {\mathcal R}_1(-\sigma)
\cale (t,\sigma,{\mathcal R}_1(\sigma)\yvec+
\frac{1}{\omega_0}{\mathcal R}_2(\sigma)\uvec)\, d\sigma
\cdot\Nabla_u \, G=0,
\\
\ds
G(t=0) = \frac{\omega_0}{2\pi}f_0,
\end{array}
\right.
\end{equation}
if $\omega_1/\omega_0 \notin \QQ$ or if $H_1=0$ and
\begin{equation}
\label{eqGder}
\left \{
\begin{array}{ll}\ds
\frac{\partial G}{\partial t} + \frac{1}{\omega_0}
\int_0^{2\pi/\omega_0} {\mathcal R}_2(-\sigma)
\Big[\cale (t,\sigma,{\mathcal R}_1(\sigma)\yvec+
\frac{1}{\omega_0}{\mathcal R}_2(\sigma)\uvec)
\\ \ds ~~~~~~~~~~~~~~~~~~~
+\frac{\omega_0}{2\pi}H_1(\omega_1 \sigma) \kkvec \big( 
{\mathcal R}_1(\sigma)\yvec+\frac{1}{\omega_0}{\mathcal R}_2(\sigma)\uvec\big)
\Big]\, d\sigma
\cdot\Nabla_y \, G
\\ \ds~~~~~~~~
+\int_0^{2\pi/\omega_0} {\mathcal R}_1(-\sigma)
\Big[\cale (t,\sigma,{\mathcal R}_1(\sigma)\yvec+
\frac{1}{\omega_0}{\mathcal R}_2(\sigma)\uvec)
\\ \ds ~~~~~~~~~~~~~~~~~~~
+\frac{\omega_0}{2\pi}H_1(\omega_1 \sigma) \kkvec \big( 
{\mathcal R}_1(\sigma)\yvec+\frac{1}{\omega_0}{\mathcal R}_2(\sigma)\uvec\big)
\Big]\, d\sigma
\cdot\Nabla_u \, G=0,
\\
\ds
G(t=0) = \frac{\omega_0}{2\pi}f_0,
\end{array}
\right.
\end{equation}
if $\omega_1/\omega_0 \in \NN$, where $\cale =\cale(t,\tau,\xvec)$ is solution to
\begin{equation}
\label{eqE}
\left \{
\begin{array}{ll}
\cale =-\nabla\Psi,
\\
\\
\ds
-\Delta \Psi= \int_{\R_v^2} G \big( t, {\mathcal R}_1(-\tau)\xvec +
\frac{1}{\omega_0}{\mathcal R}_2(-\tau)\vvec 
\, , \, - \omega_0{\mathcal R}_2(-\tau)\xvec+
{\mathcal R}_1(-\tau)\vvec\big) \, d\vvec.
\end{array}
\right.
\end{equation}
\end{theorem}
\proof
The deduction of this theorem follows directly from theorem \ref{thcali}
using the expression (\ref{exptauM}) of $e^{\tau M}$ once the two-scale
convergence in $A^\e$ given by (\ref{Aeps}) and in the Poisson equation 
(\ref{VP-paraxial}c) is achieved.

As a direct consequence of estimates (\ref{e1}) and  (\ref{e2}) and of
the regularization properties of the Laplace operator, we deduce
that $\evec_\e$ is bounded in  $ L^\infty([0,T];W^{1,3/2}(\R^2_x))$,
then, extracting a subsequence, $\evec_\e$ two-scale converges to 
$\cale= -\nabla\Psi \in L^\infty([0,T]\times[0,\frac{2\pi}{\omega_0}];W^{1,3/2}(\R^2_x))$.
\newline
In order to pass to the two-scale limit in the Poisson equation, we multiply
it by a test function $\tfct (t,\frac t\e, \xvec)$ such that 
$\tau \mapsto \tfct (t,\tau,\xvec)$ is $2\pi/\omega_0-$periodic, to give
$$
\int_0^T \hspace{-5pt} \int_{\RR^2_x}\nabla\Phi_\e(t,\xvec) \cdot 
\nabla \tfct (t,\frac t\e, \xvec) \, d\xvec dt  
= \int_0^T \hspace{-5pt} \int_{\Omega} f_\e(t,\xvec,\vvec) \,
\tfct (t,\frac t\e, \xvec) \,  d\xvec d\vvec dt,
$$
in which we pass to the limit to obtain 
\begin{multline*}
\int_0^{\frac{2\pi}{\omega_0}} \hspace{-5pt}\int_0^T \hspace{-5pt}\int_{\RR^2_x}  
\hspace{-3pt} 
\nabla\Psi (t,\xvec) \cdot 
\nabla \tfct (t,\tau, \xvec) \, d\xvec dt d\tau
=
\int_0^{\frac{2\pi}{\omega_0}}\hspace{-5pt} \int_0^T 
\hspace{-5pt}\int_{\Omega} 
F(t,\tau,\xvec,\vvec)  \tfct (t,\tau, \xvec) \, d\xvec d\vvec  dt d\tau
\\=
\int_0^{\frac{2\pi}{\omega_0}}\hspace{-5pt} \int_0^T 
\hspace{-5pt}\int_{\Omega} 
G \left( t,{\mathcal R}_1(-\tau)\xvec +
{\mathcal R}_2(-\tau)\vvec /\omega_0
\, , \, - \omega_0{\mathcal R}_2(-\tau)\xvec+
{\mathcal R}_1(-\tau)\vvec \right ) 
\tfct (t,\tau, \xvec) \,
d\xvec d\vvec dt d\tau ,
\end{multline*}
which is the weak formulation of (\ref{eqE}).

On the other hand, if $\omega_1/\omega_0 \notin \QQ$, 
for any regular function $\tfct$ which is 
$2\pi/\omega_0-$periodic in $\tau$,
the product
$$
\int_0^T \hspace{-5pt} \int_{\RR^2_x}H_1(\omega_1 \frac t\e ) \kkvec \xvec 
\cdot \tfct (t,\frac t\e, \xvec) \, d\xvec  dt\rightarrow 0,
$$
since $\tau \mapsto H_1(\omega_1 \tau)$ is $2\pi/\omega_1-$periodic
with mean value $0$.
This yields equation (\ref{eqG}).
\newline 
If $\omega_1/\omega_0 \in \NN$, then 
$\tau \mapsto H_1(\omega_1 \tau)$ is $2\pi/\omega_0-$periodic
and
$$
\int_0^T \hspace{-5pt} \int_{\RR^2_x} H_1(\omega_1 \frac t\e ) \kkvec \xvec 
\cdot\tfct (t,\frac t\e, \xvec) \, d\xvec  dt\rightarrow 
\int_0^{\frac{2\pi}{\omega_0}}\hspace{-5pt} \int_0^T \hspace{-5pt} \int_{\RR^2_x}
H_1(\omega_1 \tau) \kkvec \xvec 
\cdot \tfct (t,\tau, \xvec) \, d\xvec  dt d\tau,
$$
yielding (\ref{eqGder}) and then ending the proof.
\qed

\begin{remark} We can bring back the case 
$\omega_1/\omega_0 \in \QQ\setminus\NN$ to the case 
$\omega_1/\omega_0 \in \NN$ by finding two integers  $k$ and $l$
such that $\omega'_0 =\omega_0 /k$ and $\omega'_1 =\omega_1 /l$
satisfies $\omega'_1/\omega'_0 \in \NN$ and replacing in $(\ref{eqGder})$
$\ds \int_0^{2\pi/\omega_0 \;}d\sigma$ by 
$\ds \int_0^{2\pi/\omega'_0  \;}d\sigma$.
\end{remark}

\bigskip
We now turn to homogenizing system (\ref{2}). For simplicity,
we assume here that $\omega_0=H_0=1$. 
We also assume that $f_0$  satisfies
\begin{equation}
f_0\geq0,f_0 \in (L^1\cap L^p)(\R^2;rdrdv_r) \text{ for } p\geq 2
\text{ and }
\int_{\R^2}(r^2+v_r^2)f_0 \, rdrdv_r < + \infty.
\end{equation}
The considered system enters the framework presented in the beginning of the 
section with $\xvec$ replaced by $(r,v_r)\in\R^2$ with
\begin{equation}\label{Zut1} 
\avec = 
\begin{pmatrix}
0 \\ \evecr_\e(t, r)+H_1\hspace{-2pt}\left (\omega_1\frac{t}\e\right )r
\end{pmatrix},~~~~~
\lvec=\begin{pmatrix} v_r \\ -r\end{pmatrix},
\end{equation}
\begin{equation}\label{Zut2}
M =\begin{pmatrix} 0&1\\ -1 &0 \end{pmatrix},~~~~~~
e^{\tau M}=
\left( \begin{array}{cc}
\cos(\tau) & \sin(\tau) 
\\
-\sin(\tau) & \cos(\tau)
\end{array} \right),
\end{equation}
and the solution $f_\e$ satisfies the needed estimates to deduce the 
following Theorem.

\begin{theorem}
\label{main2}
Under the assumptions above, 
extracting from a sequence of solution $(f_\e,\evecr_\e)$ to
$(\ref{2})$ a subsequence, we deduce that  $f_\e$ two-scale converges to $F
\in L^\infty([0,T]\times[0,2\pi];L^2(\R^2;rdrdv_r))$ 
and $\evecr_\e$ two-scale converges to 
$\cale \in L^\infty([0,T]\times[0,2\pi];W^{1,3/2}(\R;rdr))$.
\newline Moreover, there exists a function 
$G=G(t,q,u_r)\in L^\infty([0,T];L^2(\R^2;qdqdu_r))$ such that
\begin{equation}
\label{Zut4}
F(t,\tau,r,v_r)= G ( t,\cos(\tau)r -\sin(\tau)v_r,\sin(\tau)r+\cos(\tau)v_r),
\end{equation}
and $G$ is solution to
\begin{equation}
\label{Zut5}
\left \{
\begin{array}{ll}\ds
\frac{\partial G}{\partial t} + 
\int_0^{2\pi} -\sin(\sigma)\,
\caler (t,\sigma,\cos(\sigma)q+\sin(\sigma)u_r)\, d\sigma
\; \frac{\partial G}{\partial q} 
\\ \ds~~~~~~~~~
+\int_0^{2\pi} \cos(\sigma)\,
\caler (t,\sigma,\cos(\sigma)q+\sin(\sigma)u_r)\, d\sigma
\; \frac{\partial G}{\partial u_r}  =0,
\\
\ds
G(t=0) = \frac{1}{2\pi}f_0,
\end{array}
\right.
\end{equation}
if $\omega_1/\omega_0 \notin \QQ$ or if $H_1=0$ and 
\begin{equation}
\label{Zut5bis}
\left \{
\begin{array}{ll}\ds
\frac{\partial G}{\partial t} + 
\int_0^{2\pi} -\sin(\sigma)\bigg(
\caler (t,\sigma,\cos(\sigma)q+\sin(\sigma)u_r)
\\ \ds ~~~~~~~~~~~~~~~~~~~~~~~~~
+ \frac{1}{2\pi} H_1(\omega_1\sigma)\big(\cos(\sigma)q+\sin(\sigma)u_r \big)
\bigg)\, d\sigma
\; \frac{\partial G}{\partial q} 
\\ \ds~~~~~~~~~
+\int_0^{2\pi} \cos(\sigma)\bigg(
\caler (t,\sigma,\cos(\sigma)q+\sin(\sigma)u_r)
\\ \ds ~~~~~~~~~~~~~~~~~~~~~~~~~
+ \frac{1}{2\pi} H_1(\omega_1\sigma)\big(\cos(\sigma)q+\sin(\sigma)u_r \big)
\bigg)\, d\sigma
\; \frac{\partial G}{\partial u_r}  =0,
\\
\ds
G(t=0) = \frac{1}{2\pi}f_0,
\end{array}
\right.
\end{equation}
if $\omega_1/\omega_0 \in \NN$
where $\caler=\caler(t,\tau,r,v_r)$ is given by
\begin{equation}
\label{Zut6}
\frac 1r \frac{\partial (r\caler)}{\partial r} 
= {\Upsilon}(t,\tau, r)=\int_\R G \big( t, \cos(\tau)r -\sin(\tau)v_r
,\sin(\tau)r +\cos(\tau)v_r )\, dv_r \, .
\end{equation}

\end{theorem}

\section{The two scale PIC solver}
\label{numerics}

In this section, we develop a two-scale PIC-method, tailored to approximate $f_\e$ very efficiently on long time scales, in the simplified case of axisymmetric beams.
The strategy consists in computing the solution $F$ to 
(\ref{Zut4})-(\ref{Zut5})-(\ref{Zut6}) or 
(\ref{Zut4})-(\ref{Zut5bis})-(\ref{Zut6}) and then to approach the solution
$f_\e$ of (\ref{2}) by $F(t,t/\e,\xvec,\vvec)$.
The advantage of proceeding in such a way is that the solution of 
(\ref{Zut4})-(\ref{Zut5})-(\ref{Zut6}) or 
(\ref{Zut4})-(\ref{Zut5bis})-(\ref{Zut6}) does not contain $1/\eps-$frequency
oscillations. As a consequence, a much larger time step may be used in the numerical method.
First, we present the implemented algorithm. Then, we compare
the solution obtained with our method with the solution 
$f_\e$ directly computed from (\ref{2}).
Finally, we compare the performances of both methods.

\subsection{Description of the numerical method}
 
The Particle In Cell (PIC) method for the Vlasov-Poisson (or Vlasov-Maxwell) equations consists in approximating the distribution function defined in phase space by a meshless particle method and the electric field on a grid of the physical space only. In addition, adequate interpolation and charge deposition methods are used to transfer the needed quantities between grid and particles.
 
For the two-scale method we want to develop, the Vlasov-Poisson equation are replaced by equations (\ref{Zut4})-(\ref{Zut5})-(\ref{Zut6}) or 
(\ref{Zut4})-(\ref{Zut5bis})-(\ref{Zut6}).

The key point of the algorithm is the computation of function $G$ solution to
 (\ref{Zut5})-(\ref{Zut6}) 
or (\ref{Zut5bis})-(\ref{Zut6}) at time $t_{l+1}=t_l+\Delta t$
knowing it at time $t_l$. 
We give the details of  the method only in the case of system (\ref{Zut5})-(\ref{Zut6}),
knowing that in the case of system (\ref{Zut5bis})-(\ref{Zut6}) ad-hoc
terms need to be managed.

As usual in a PIC-method, $G$ is 
approximated by the following Dirac mass sum
\begin{equation} G_N(q,u,t)=\sum_{k=1}^N w_k\delta(q-Q_k(t))\delta(u-U_k(t)),\end{equation}
where $(Q_k(t),U_k(t))$ is the position in phase space of macro-particle $k$ which moves
along a characteristic curve of the first order PDE (\ref{Zut5}).
Hence the job is reduced to compute the macro-particle positions
$(Q_k^{l+1},U_k^{l+1})$ at time $t_{l+1}=t_l+\Delta t$ 
from their positions $(Q_k^{l},U_k^{l})$ at time $t_l$,
knowing they are solutions to
\begin{align}
&\frac{d Q_k}{dt}=-\int_0^{2\pi}\sin(\sigma)\,
\caler (t,\sigma,\cos(\sigma)Q_k+\sin(\sigma)U_k)\, d\sigma,
~~&Q_k(t_l)=Q_k^{l},
\label{Zut10}
\\
&\frac{d U_k}{dt}=\int_0^{2\pi}\cos(\sigma)\,
\caler (t,\sigma,\cos(\sigma)Q_k+\sin(\sigma)U_k)\, d\sigma,
~~&U_k(t_l)=U_k^{l}.
\label{Zut10-1}
\end{align}
These characteristic equations are a lot more complex to deal with than those of the standard Vlasov-Poisson system which read $\frac{dX}{dt}=V$,   $\frac{dV}{dt}=E(X,t)$.
Indeed, in our case $Q$ and $U$ are coupled in each equation and an integral term needs to be approximated.

Because of the form of the right 
hand side in (\ref{Zut10}) and (\ref{Zut10-1})
all along the algorithm, we need to compute
values of the two-scale electric field $\caler$
generated by a given macro-particle distribution
$(Q_k,U_k)_{k=1,..,N}$. 
The tedious step while computing  $\caler$ in grid 
point $(q_i)_{i=1,..,A}$ is the computation of the right
hand side ${\Upsilon}(t,\sigma_m,q_i)$ of (\ref{Zut6}).
Indeed, ${\Upsilon}(t,\sigma_m,q_i)$ involves the integral of the
particle distribution on the oblique line which is the range 
of the vertical line $[q=q_i]$ by the rotation 
$e^{-\sigma_m M}$ defined by (\ref{Zut2}).
Hence we have to apply this rotation to each line 
$[q=q_i]$ and to project the particles 
$(Q_k^{},U_k^{})_{k=1,..,N}$ on the resulting oblique lines.
Summing then the projection result on the oblique line associated 
with $q_i$ yields the value of $\Upsilon(t_l,\sigma_m,q_i)$.
Another way to obtain this value consists in applying the rotation
$e^{\sigma_m M}$ to the particles, then projecting this rotation result 
on lines $[q=q_i]$ and summing.
Once $\Upsilon(t_l,\sigma_m,q_i)$  
is known in each $q_i$, the computation of the $\caler (t_l,\sigma_m,q_i)$
are straightforward using any classical Poisson numerical solver.

The first step of the  computation of $(Q_k^{l+1},U_k^{l+1})$ consists 
in replacing the integrals above by $p-$node quadrature 
formula. As we approximate the integral of a periodic function over one period, the trapezoidal rule is optimal and will yield very accurate results for as few quadrature points as are needed to resolve the oscillations of the function.

Then, the equations for $(Q_k,U_k)$ become
\begin{gather}
\label{Zut11}
\frac{d Q_k}{dt}=-\sum_{m=1}^{p}\gamma_m\sin(\sigma_m)\,
\caler (t,\sigma_m,\cos(\sigma_m)Q_k+\sin(\sigma_m)U_k),
~~Q_k(t_l)=Q_k^{l},
\end{gather}
\begin{gather}
\label{Zut11bis}
\frac{d U_k}{dt}=\sum_{m=1}^{p}\gamma_m\cos(\sigma_m)\,
\caler (t,\sigma_m,\cos(\sigma_m)Q_k+\sin(\sigma_m)U_k),
~~U_k(t_l)=U_k^{l}.
\end{gather}

Then, we solve (\ref{Zut11})-(\ref{Zut11bis}) using the classical
Runge-Kutta method:
\begin{equation}
\begin{array}{c|cccc}
0&&&&\\
1/2&1/2&&&\\
1/2&0&1/2&&\\
1&0&0&1&\\
\hline
&1/6&1/3&1/3&1/6
\end{array}~~~~
\end{equation}
which gives the following scheme when applied to the 
computation of the approximation $y^{l+1}$ of the
value of $y$ solution to $dy/dt=K(t,y)$ at time $t_l+\Delta t$
knowing its approximation $y^{l}$ at time $t_l$:
\begin{equation}
\begin{aligned}
& t_{l,1}=t_l, ~~  y^{l,1}= y^{l} 
\\ 
& t_{l,2}=t_l+\frac{\Delta t}{2}, ~~  y^{l,2}= y^{l} + \frac 12 I^1, 
\text{ with } I^1= \Delta t \, K(t_{l,1},y^{l,1}),
\\
& t_{l,3}=t_l+\frac{\Delta t}{2}, ~~  y^{l,3}= y^{l} + \frac 12 I^2, 
\text{ with } I^2= \Delta t \, K(t_{l,2},y^{l,2}),
\\
& t_{l,4}=t_l+\Delta t, ~~  y^{l,4}= y^{l} + I^3, 
\text{ with } I^3= \Delta t  \, K(t_{l,3},y^{l,3}),
\\
\\
&y^{l+1} =y^{l}+ \frac 16 I^1 + \frac 13 I^2 + \frac 13 I^3 + \frac 16 I^4 
\text{ with } I^4= \Delta t  \, K(t_{l,4},y^{l,4}).
\end{aligned}
\end{equation}
Applying this scheme to our problem consists in replacing in the 
formula above $y$ by $(Q_k,U_k)$ and computing $K$ using the
result of a Poisson solver. 
In other words, we have to compute $Q_k^{l,2}$ as follows:
\begin{align}
&Q_k^{l,2}=
Q_k^{l} +\frac{1}{2}I^1, \text{ with } \nonumber
\\
&I^1= \Delta t\Big(
-\sum_{m=1}^{p}\gamma_m\sin(\sigma_m)\,
\caler (t_l,\sigma_m,\cos(\sigma_m)Q_k^{l}+\sin(\sigma_m)U_k^{l})
                           \Big).
\end{align}
and something similar for $U_k^{l,2}$.

In order to achieve this, we compute the value of $\caler$
in $(t_l,\sigma_m,\cos(\sigma_m)Q_k^{l}+\sin(\sigma_m)U_k^{l})$ 
by interpolating the value of $\caler(t_l,\sigma_m,q_i)$ known
on the grid $(q_i)_{i=1..,A}$ as soon as it has been computed
solving the Poisson equation  (\ref{Zut6}) associated with the 
particle distribution $(Q_k^{l},U_k^{l})$ by the procedure
described above.

The following step of the Runge-Kutta method consists in computing
$Q_k^{l,3}$ defined by
\begin{align}
&Q_k^{l,3}=
Q_k^{l} +\frac{1}{2} I^2, \text{ with } \nonumber
\\
&I^2= \Delta t\Big(
-\sum_{m=1}^{p}\gamma_m\sin(\sigma_m)\,
\caler^2 (t_l+\frac{\Delta t}{2},\sigma_m,
         \cos(\sigma_m)Q_k^{l,2}+\sin(\sigma_m)U_k^{l,2})
                           \Big),
\end{align}
where the value of 
$\caler^2 (t_l+\frac{\Delta t}{2},\sigma_m,
\cos(\sigma_m)Q_k^{l,2}+\sin(\sigma_m)U_k^{l,2})$
is obtained by interpolation of 
$\caler^2 (t_l+\frac{\Delta t}{2},\sigma_m,q_i)$
which is computed as previously from the 
$(Q_k^{l,2},U_k^{l,2})_{k=1,..,N}$ particle distribution.

Then we compute 
\begin{align}
&Q_k^{l,4}=
Q_k^{l} + I^3, \text{ with } \nonumber
\\
&I^3 = \Delta t\Big(
-\sum_{m=1}^{p}\gamma_m\sin(\sigma_m)\,
\caler^3 (t_l+\frac{\Delta t}{2},\sigma_m,
         \cos(\sigma_m)Q_k^{l,3}+\sin(\sigma_m)U_k^{l,3})
                           \Big),
\end{align}
where $\caler^3 (t+\frac{\Delta t}{2})$ is computed from particle positions
$(Q_k^{l,3},U_k^{l,3})_{k=1,..,N}$.

Finally, $Q_k^{l+1}$ is obtained by the following formula:
\begin{align}
&Q_k^{l+1}=Q_k^{l} + \frac{1}{6} I^1 +\frac{1}{3} I^2 
+\frac{1}{3} I^3 +\frac{1}{6} I^4 , \text{ with } \nonumber
\\
&I^4 = \Delta t \Big(
-\sum_{m=1}^{p}\gamma_m\sin(\sigma_m)\,
\caler^4 (t_l+\Delta t,\sigma_m,
         \cos(\sigma_m)Q_k^{l,4}+\sin(\sigma_m)U_k^{l,4})
                 \Big).
\end{align}
where $I^1$, $I^2$ and $I^3$ are defined above
and where $\caler^4(t_l+\Delta t)$ is computed from particle positions 
$(Q_k^{l,4},U_k^{l,4})_{k=1,..,N}$.

\subsection{Validation of the two-scale solver}

\paragraph{The linear case.} First in order to check our implementation
we consider the case when the self-consistent electric field vanishes.
Then by choosing adequately the function $H_1$ we can compute an
analytical solution of the two-scale model that we can compare with our code
and with the solution given by the usual PIC solver.

\begin{figure}
\begin{center}
\begin{tabular}{ll}
\includegraphics[height=4cm]{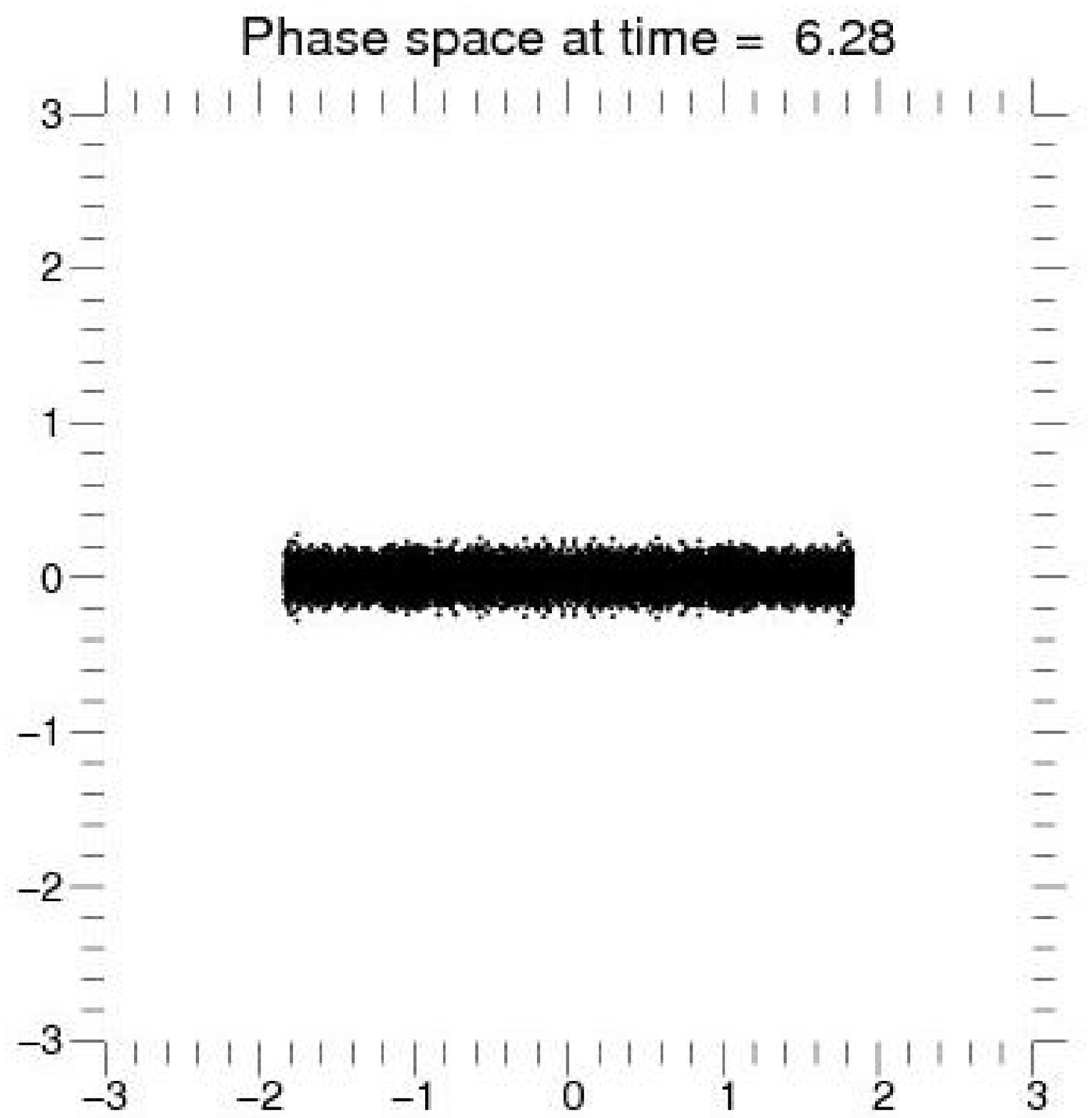}
&
\includegraphics[height=4cm]{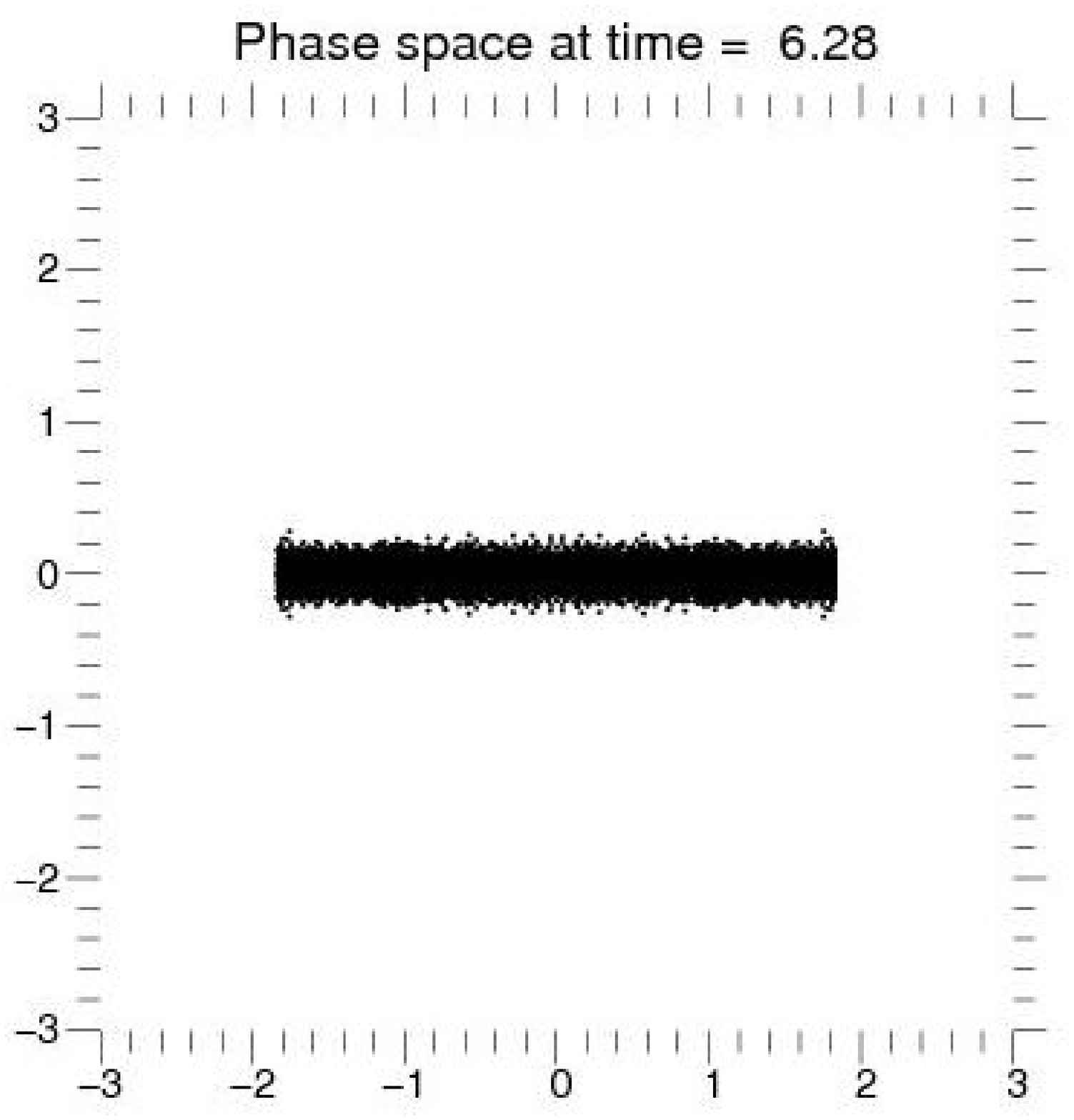}
\end{tabular}
\caption{\label{linNR} Beam simulation without self-consistent electric field in the non resonant case ($H_1(t) = \cos(4\sqrt {2}t)$) with a usual PIC method (left) and a two-scale 
PIC method (right) for $\eps = 0.01$ at time 6.28. }
 \end{center}
\end{figure}

\begin{figure}
\begin{center}
\begin{tabular}{ll}
\includegraphics[height=4cm]{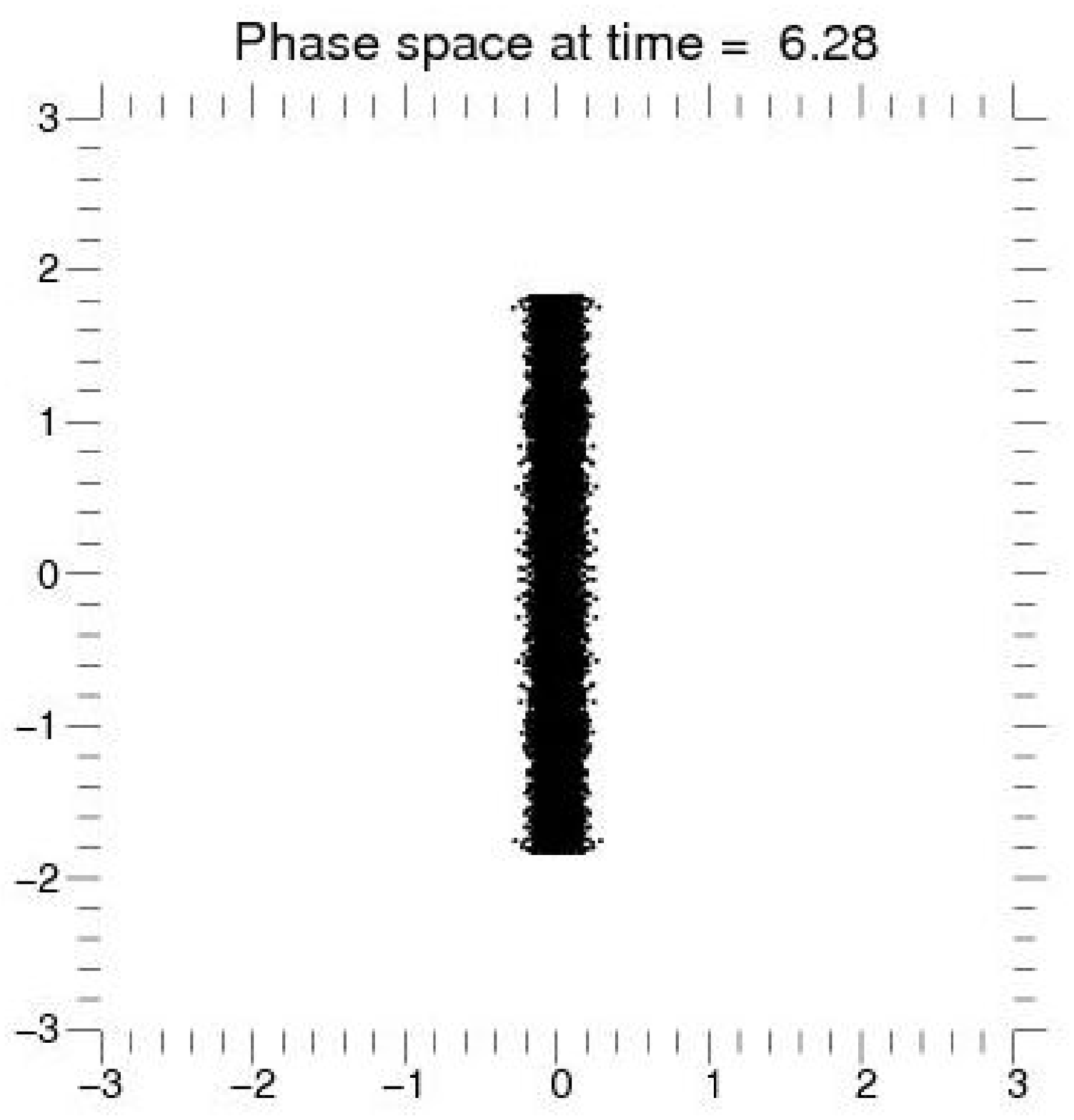}
&
\includegraphics[height=4cm]{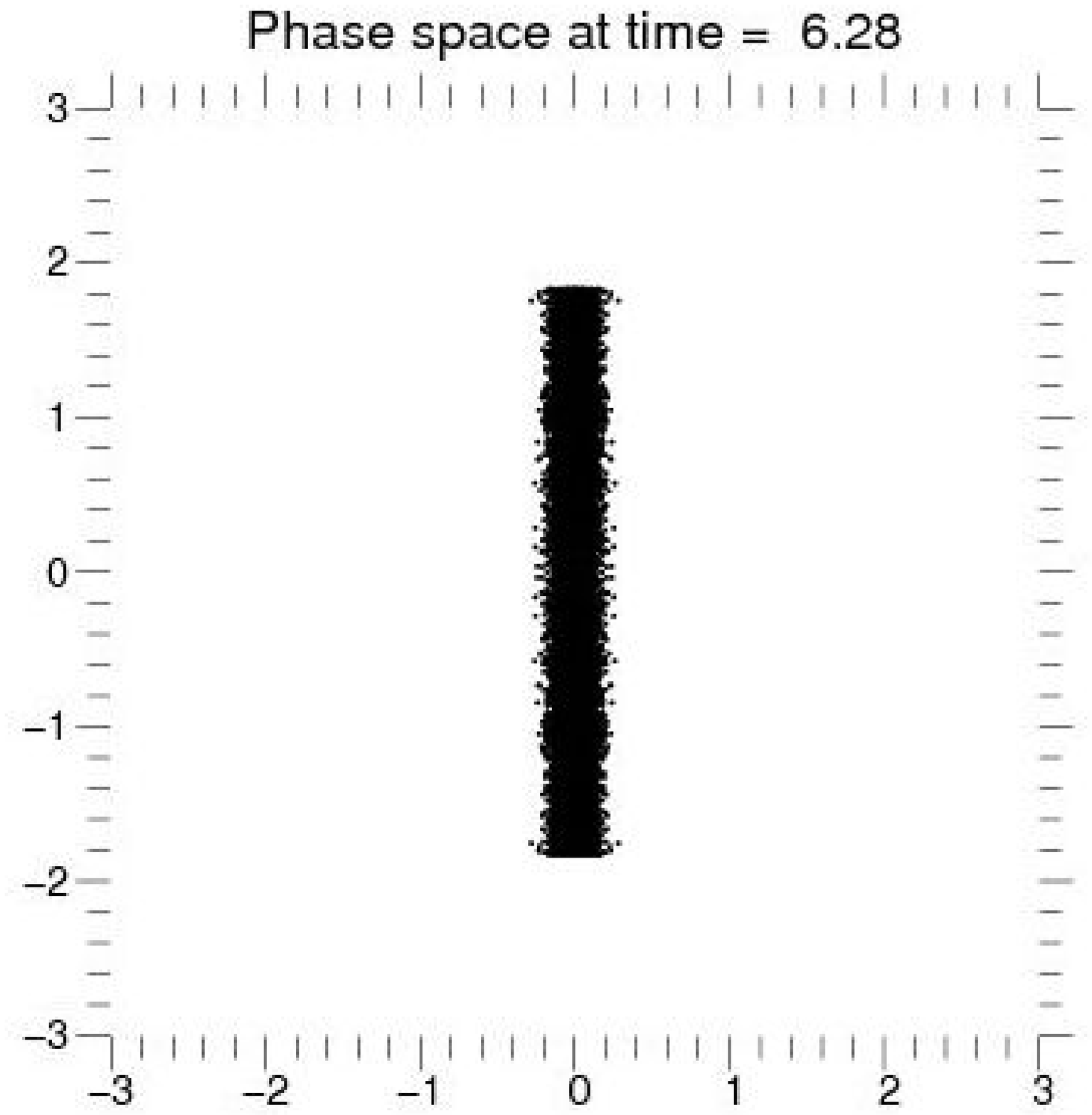}
\end{tabular}
\caption{\label{linR} Beam simulation without self-consistent electric field in the resonant case ($H_1(t)=\cos^2(2t)$) with a usual PIC method (left) and a two-scale 
PIC method (right) for $\eps = 0.01$ at time 6.28. }
 \end{center}
\end{figure}

Let us first consider the non resonant case by choosing $H_1(t) = \cos(4\sqrt {2}t)$.
In this case, as the self-consistent electric field vanishes, $G$ is stationary and the
beam will only move through the rotation transforming $G$ to $F$. We check on
figure \ref{linNR} that this is indeed the case. Let us now consider the resonant
case by choosing $H_1(t)=\cos^2(nt)$ with $n\geq 2$. Then a straightforward computation yields the equation satisfied by $G$ 
$$\frac{\partial G}{\partial t} -\frac 14 u \frac{\partial G}{\partial q} +\frac 14 q \frac{\partial G}{\partial u}.$$
The characteristics of this equation can be computed explicitly:
$$
Q(t) = Q_0\cos \frac t4 - U_0\sin \frac t4, \quad
U(t) = Q_0\sin \frac t4 + U_0\cos \frac t4, 
$$
which yields an explicit solution of $G$ and also an explicit solution of $F$ using formula
\eqref{Zut4}. Comparing this exact solution with the solution computed by our solver we can check its accuracy. Let us first mention, that our quadrature formula yields the exact result up to machine accuracy for any number of quadrature points greater or equal to seven for $n=2$ in the definition of $H_1$. Of course when the chosen $n$ is larger the function oscillates more and more quadrature points are needed. For $n=7$, for example, the exact result, up to machine accuracy, is obtained with 17 quadrature points. We can thus conclude that as long as the oscillations are resolved our quadrature rule is very accurate. We also check that our RK solver is of order 4 in $\Delta t$ as expected.
In figure  \ref{linR}, we check the global behavior for a whole beam in comparison
to the usual PIC solver and the results are very satisfying. The phase shift appearing in the two-scale limit is indeed also observed in the usual PIC code. Note that because the usual PIC code needs to resolve the fast time scale, the difference in time step between of the two solvers is of the order of $\epsilon$, so that the efficiency of the two-scale PIC solver compared to the usual PIC solver is better for small values of $\epsilon$.

\paragraph{The non linear case.} In this case, as no analytical solution is available,
we shall compare our solver with a traditional PIC solver resolving the small time scale.

Let us first consider a case of a fairly small $\eps$  on which we test our two-scale PIC method, where the two scale solution should be fairly close to the solution given
by the traditional PIC solver. 
We take $\eps=0.01$, and choose as an initial distribution 
\begin{equation}
f_0(r,v_r) = 
\frac{n_0}{\sqrt{2\pi} v_{th}} 
exp\bigg( -\frac{v_r^2}{2v_{th}^2} \bigg)
\chi_{[-0.75,0.75]}(r), 
\end{equation}
with thermal velocity $v_{th}= 0.0727518214392$ and where $\chi_{[-0.75,0.75]}(r) =1$
if $r \in  [-0.75,0.75]$ and 0 otherwise and in considering function $H_1$ as being
identically zero. This corresponds to a semi-Gaussian beam often used in accelerator physics. We consider in all our simulations beams which are willingly unmatched in order to assess to what extent our two-scale PIC solver can follow the complex structures that are developing.

\begin{figure}
\begin{center}
\hspace{-5mm}
\begin{tabular}{lll}
\includegraphics[height=4cm]{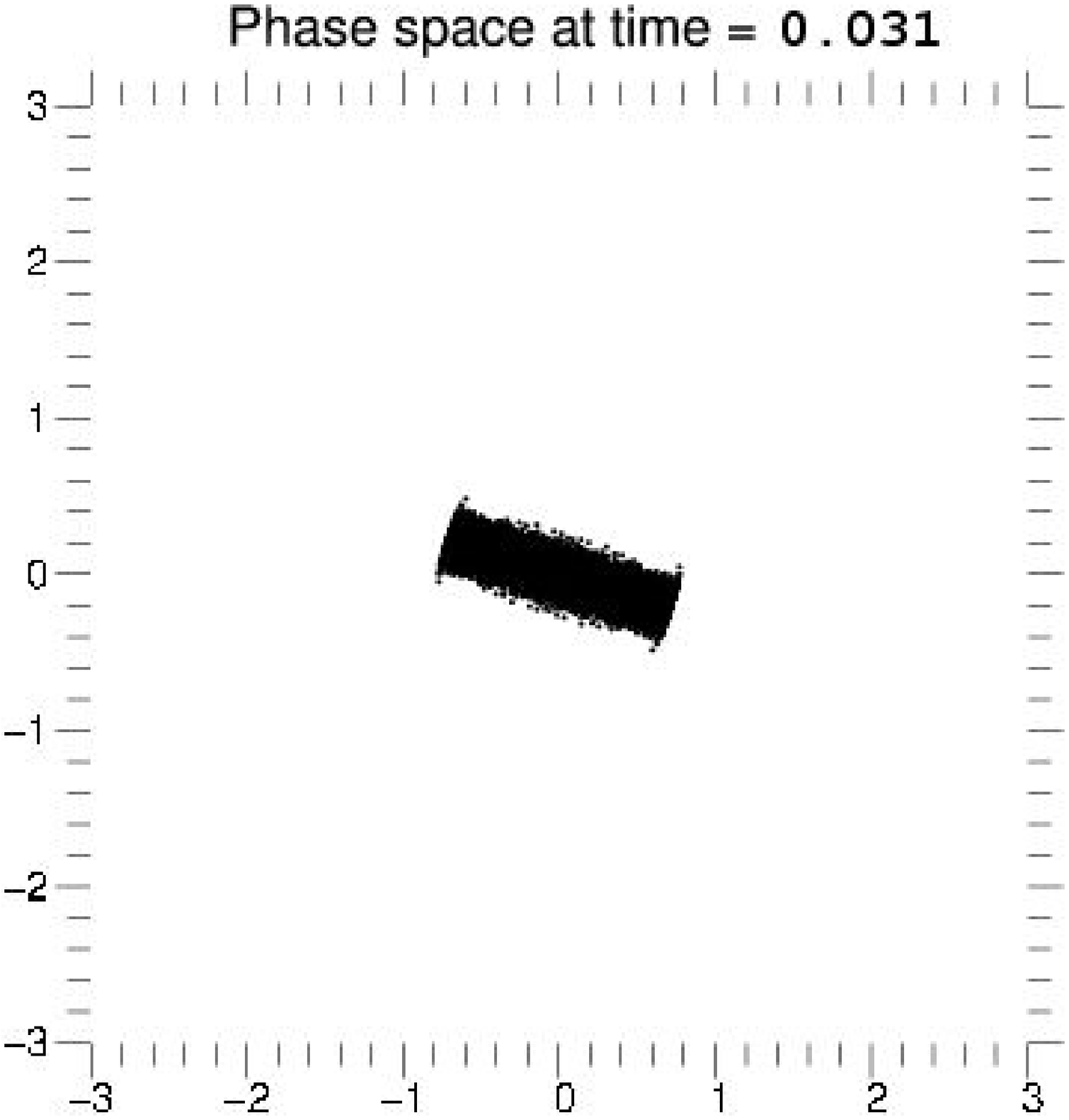}
&
\includegraphics[height=4cm]{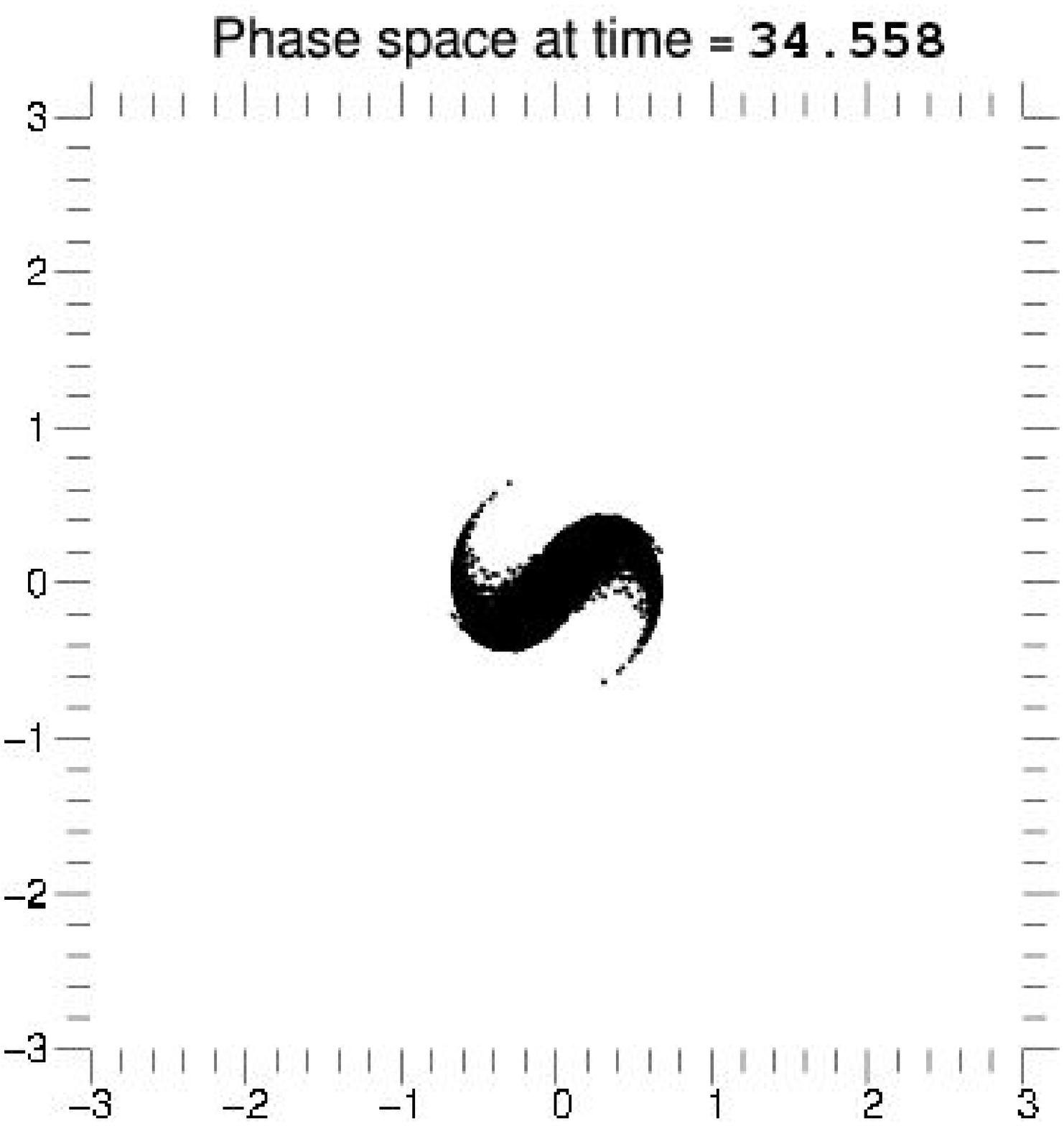}
&
\includegraphics[height=4cm]{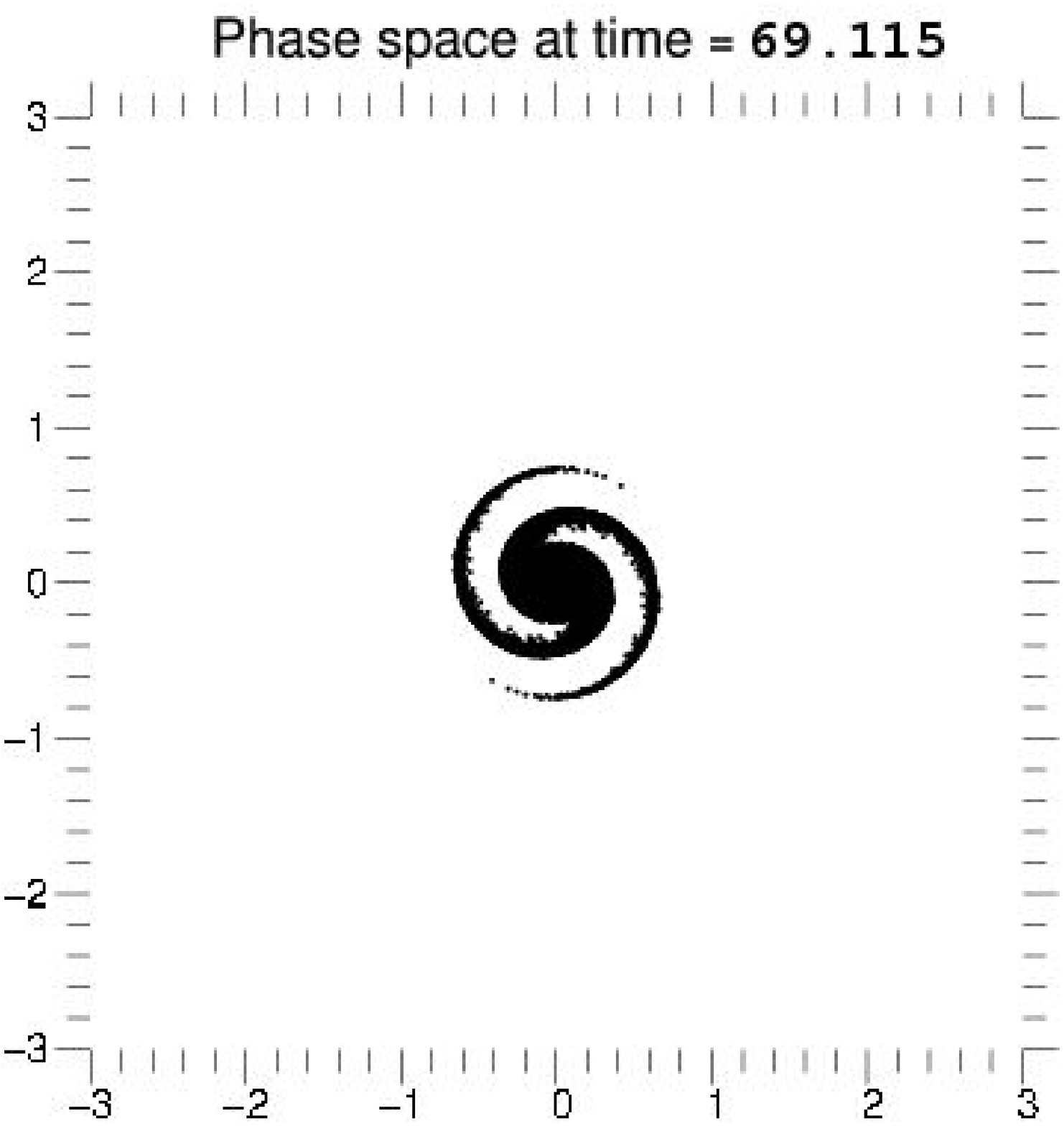}
\\ \\
\includegraphics[height=4cm]{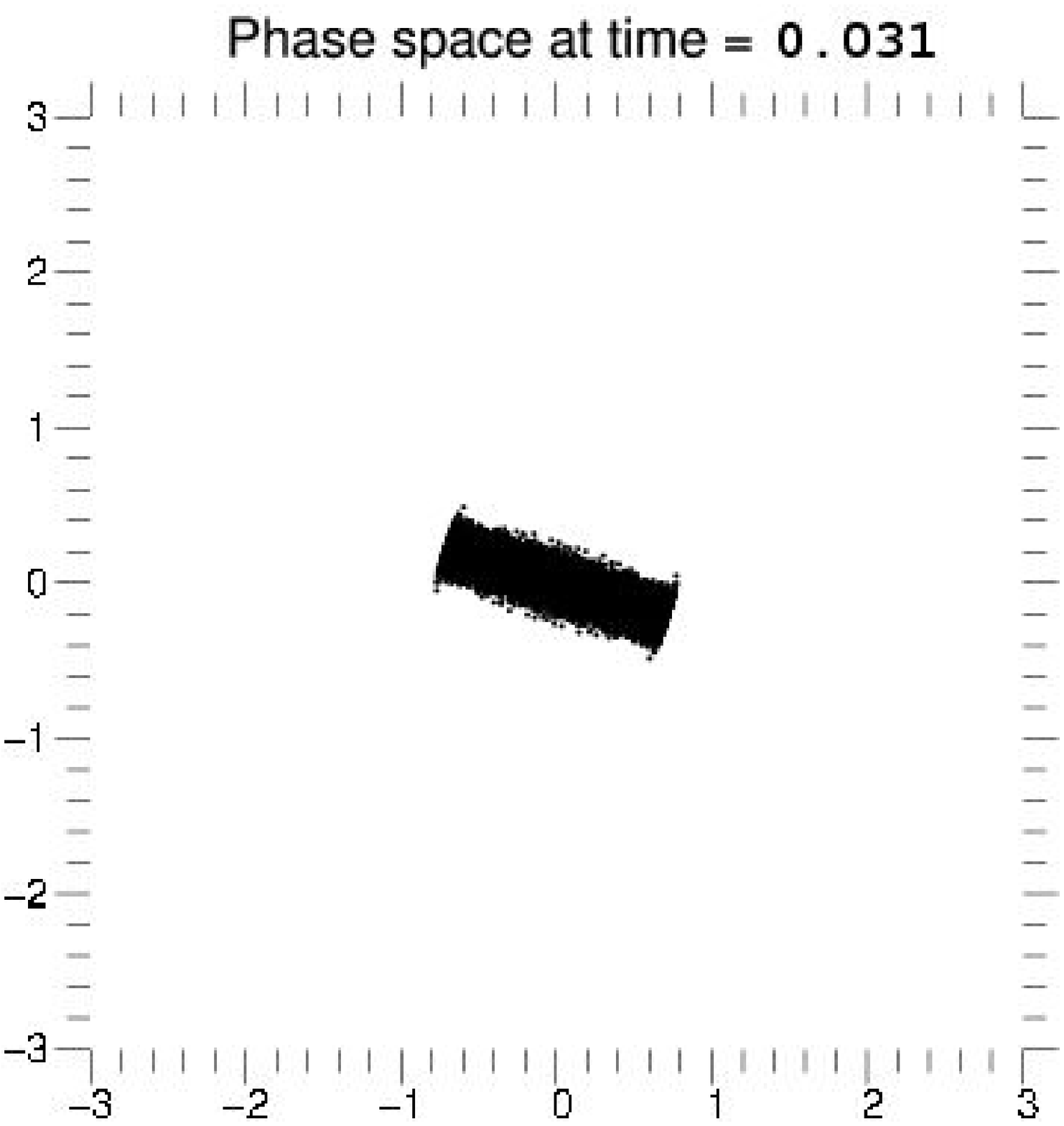}
&
\includegraphics[height=4cm]{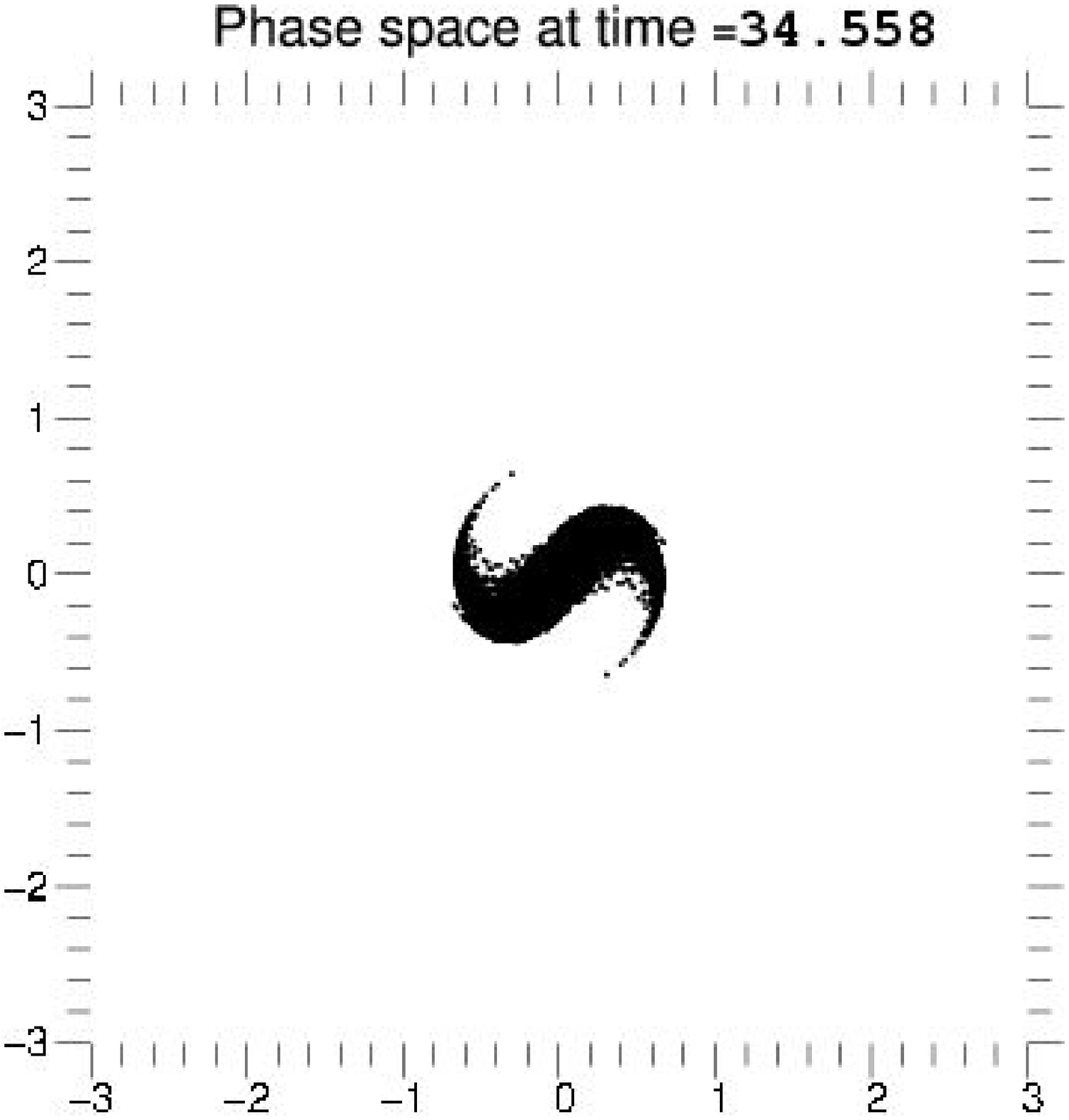}
&
\includegraphics[height=4cm]{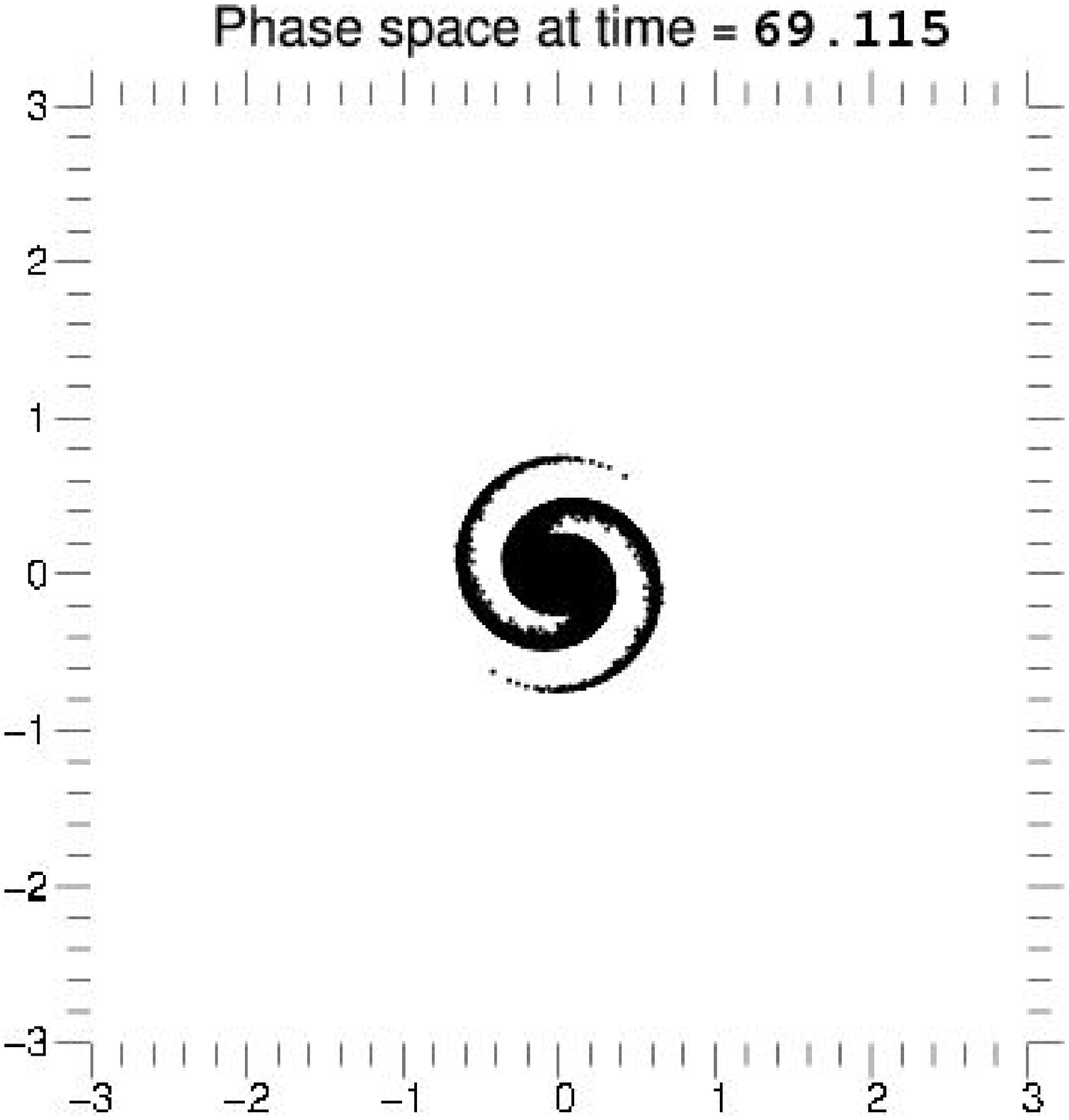}
\end{tabular}
\caption{Beam simulation with a usual PIC method and a two-scale 
PIC method for $\eps = 0.01$. Left: beam at time 0.031, center: 
beam at time 34.558, right : beam at time 69.115. 
Top : Simulation provided with the usual PIC method,
 bottom:  Simulation provided with the two-scale PIC method.}\label{fig1}
 \end{center}
\end{figure}

We use a 15-node composed trapezoidal quadrature formula.
The results are given in figure \ref{fig1} where the horizontal axis is the $r$ axis
and the vertical axis the $v_r$ axis. The top line shows the time evolution
of the beam simulated with a usual standard PIC method and the bottom line
shows the same simulation with the two-scale PIC method just built.
We can see that the two simulations coincide with a high degree of accuracy.
The time step in the two-scale PIC method is $\eps$ times larger
as in the usual PIC method.
The simulations were both made with an Intel Core 2 Duo processor (2.33 GHz)
under  Mac OS X 10.4.10 (8R2218) system. The needed CPU time for the simulation
with the usual standard PIC method is 4320.382 seconds and is 549.197 seconds  with
the two-scale PIC method.

\medskip
The three next tests are done with $\eps=0.1$ and with initial distribution
\begin{equation}
f_0(r,v_r) = 
\frac{n_0}{\sqrt{2\pi} v_{th}} 
exp\bigg( -\frac{v_r^2}{2v_{th}^2} \bigg)
\chi_{[-1.83271471003,1.83271471003]}(r), 
\end{equation}
with thermal velocity $v_{th}= 0.0727518214392$.
In the case when 
\begin{equation}
H_1(\omega_1 \tau) = \cos(\tau),
\end{equation}
the mean effect of it is zero in the sense that terms containing 
$H_1$ appearing in (\ref{Zut5bis}) are both zero.
\begin{figure}
\begin{center}
\begin{tabular}{cc}
\includegraphics[height=4.7cm]{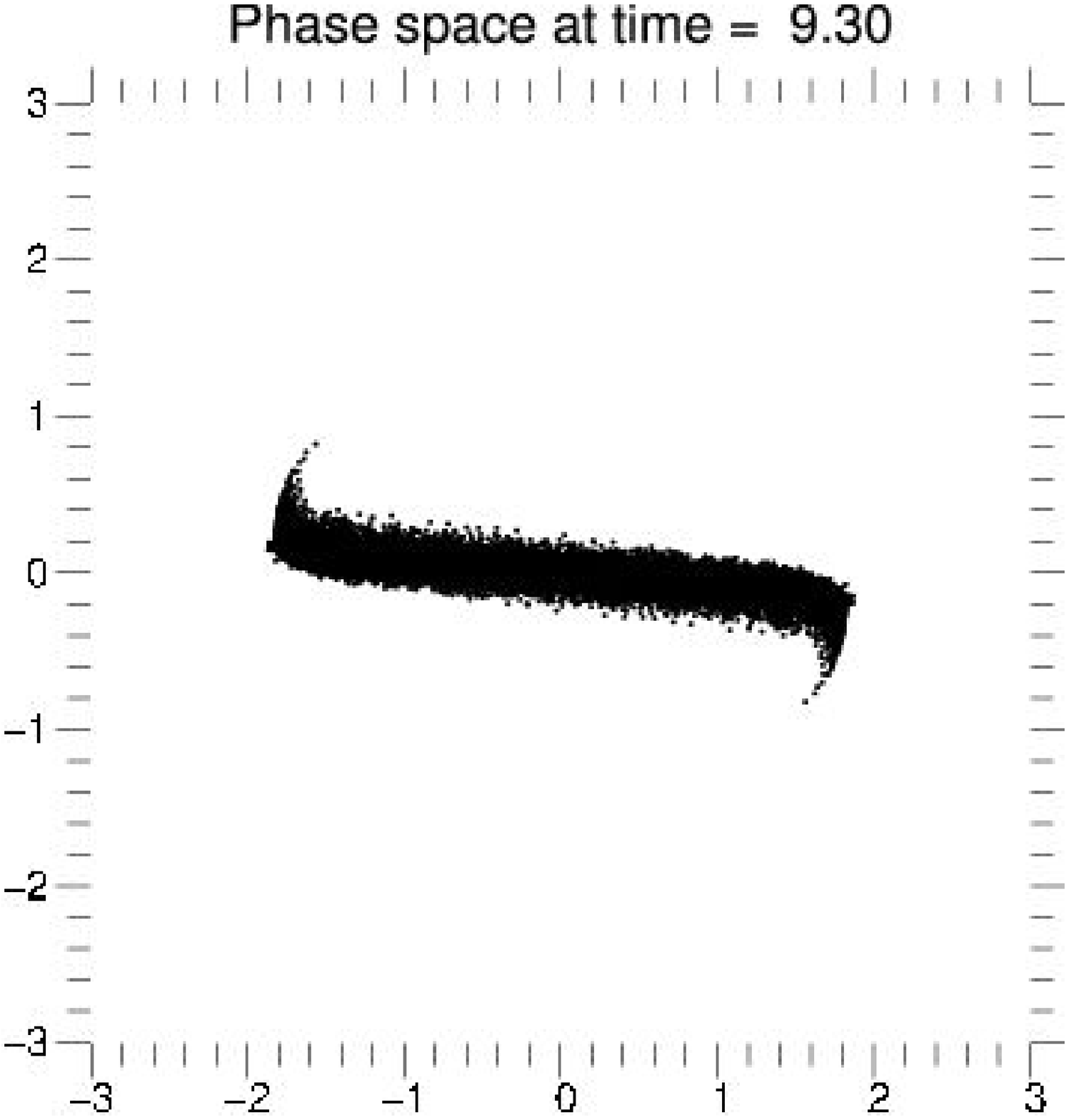}
&
\includegraphics[height=4.7cm]{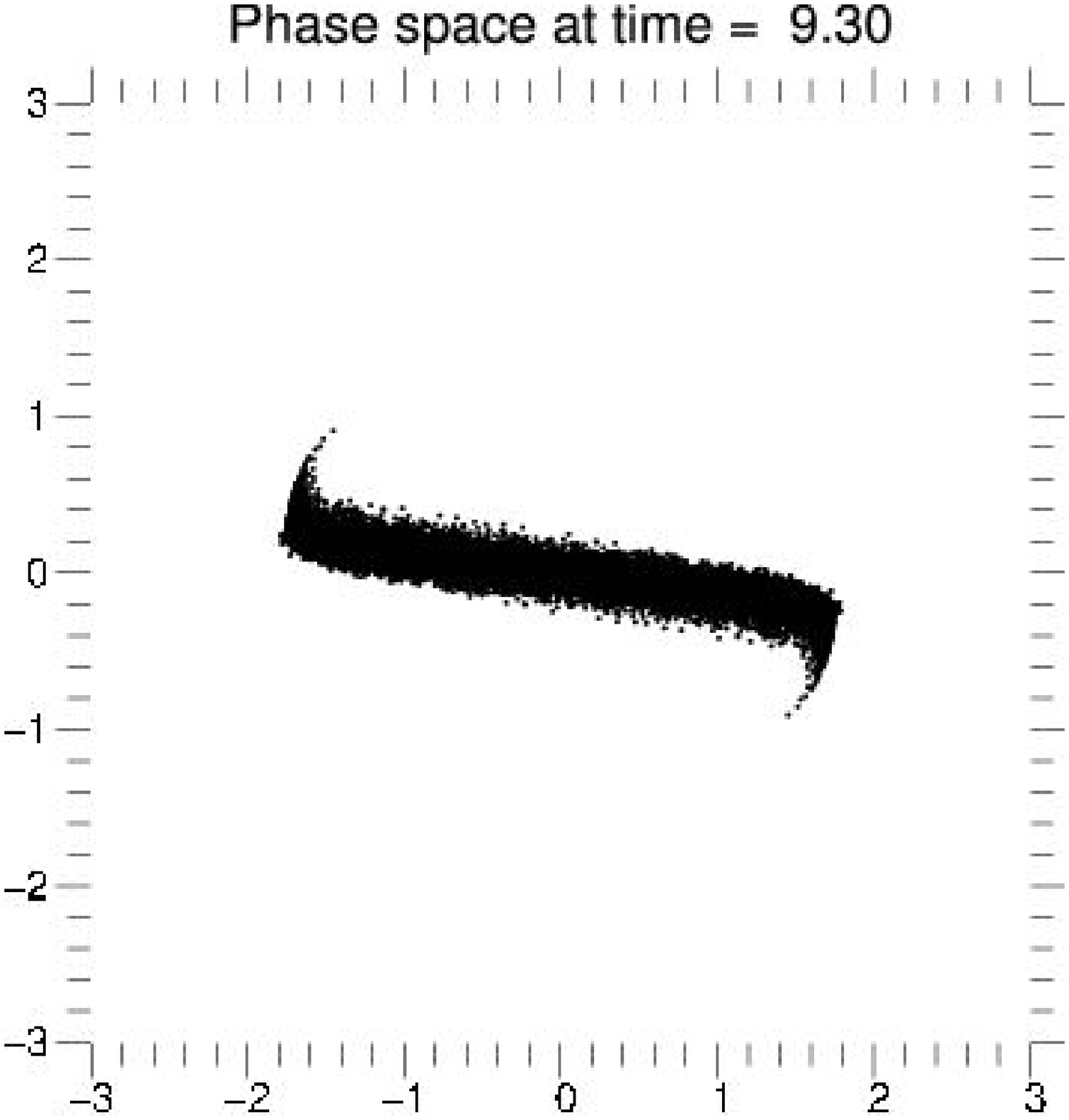}
\end{tabular}
\end{center}
\caption{Beam simulation with an external force with no mean effect at time 9.30.
Left: with the usual PIC method. Right: with the two-scale PIC method.}
\label{fig2}
\end{figure}
\begin{figure}
\begin{center}
\begin{tabular}{cc}
\includegraphics[height=4.7cm]{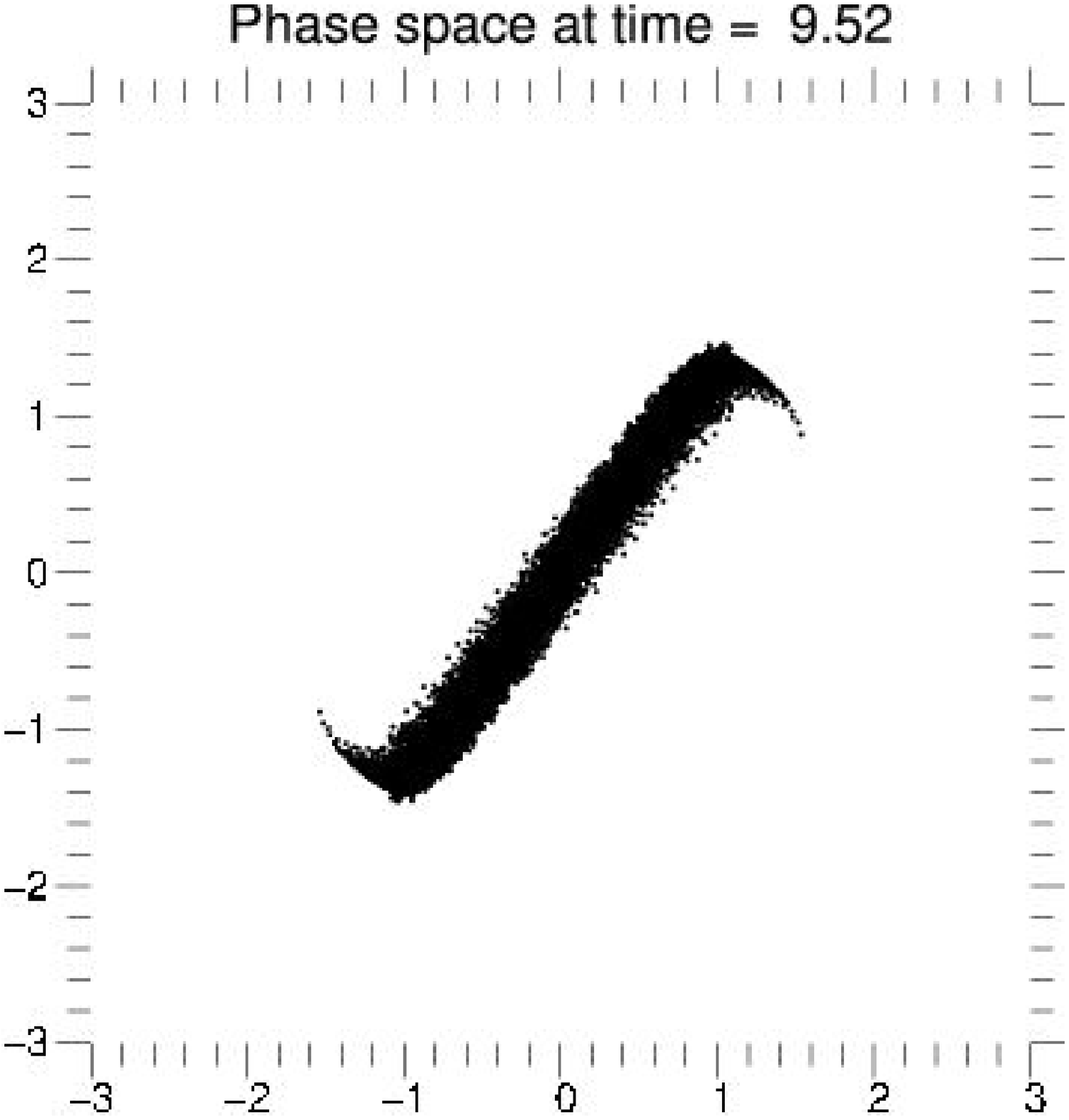}
&
\includegraphics[height=4.7cm]{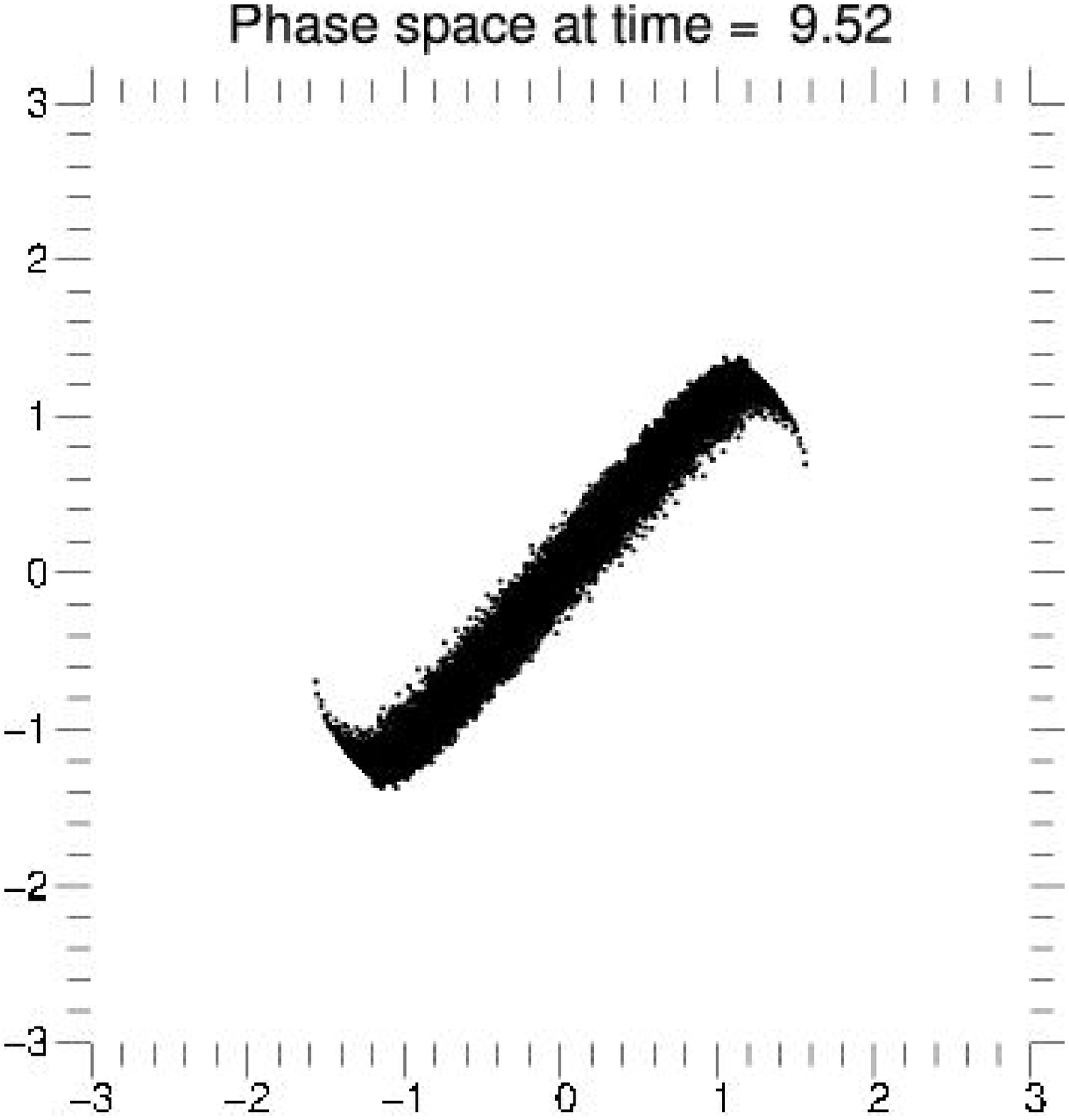}
\end{tabular}
\end{center}
\caption{Beam simulation with an external force with no mean effect at time 9.52.
Left: with the usual PIC method. Right: with the two-scale PIC method. }
\label{fig3}
\end{figure}
Numerical results concerning this case are given in figures \ref{fig2} and
\ref{fig3}.  On the left of both figures is shown the simulation result with a
usual PIC method. On the right is shown the simulation result with the 
two-scale PIC method. We can notice that if the beam configuration are very similar
at time 9.30 (figure \ref{fig2}) the two-scale PIC method is a bit in advance
at time 9.52 (figure \ref{fig3}). 
This illustrates a finite $\eps$ effect making that the 
two-scale PIC method is sometimes late and sometimes in advance with
respect to the usual PIC method. 

\begin{figure}
\hspace{-5mm}
\begin{tabular}{llll}
\includegraphics[height=4cm]{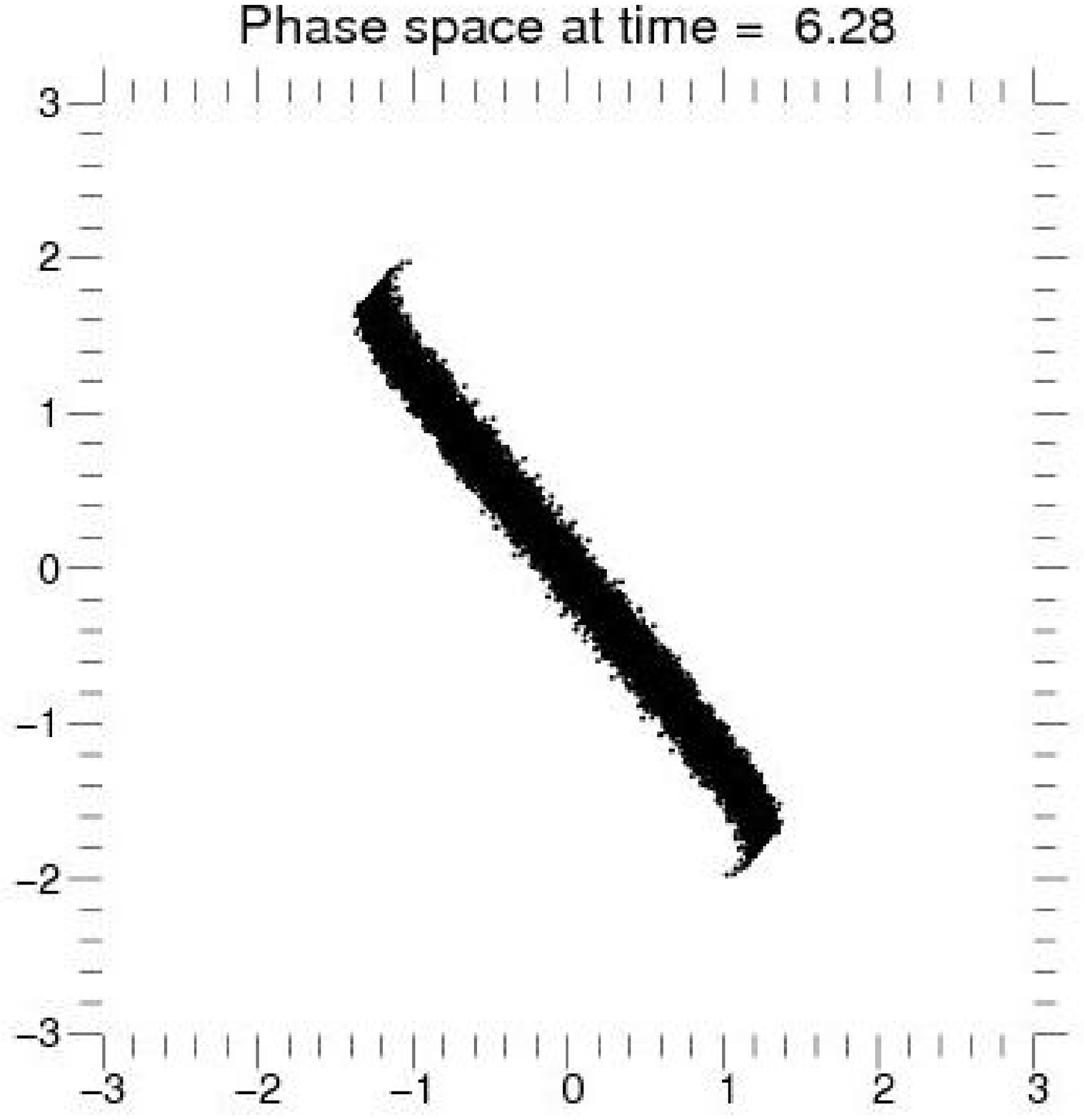}
&
\includegraphics[height=4cm]{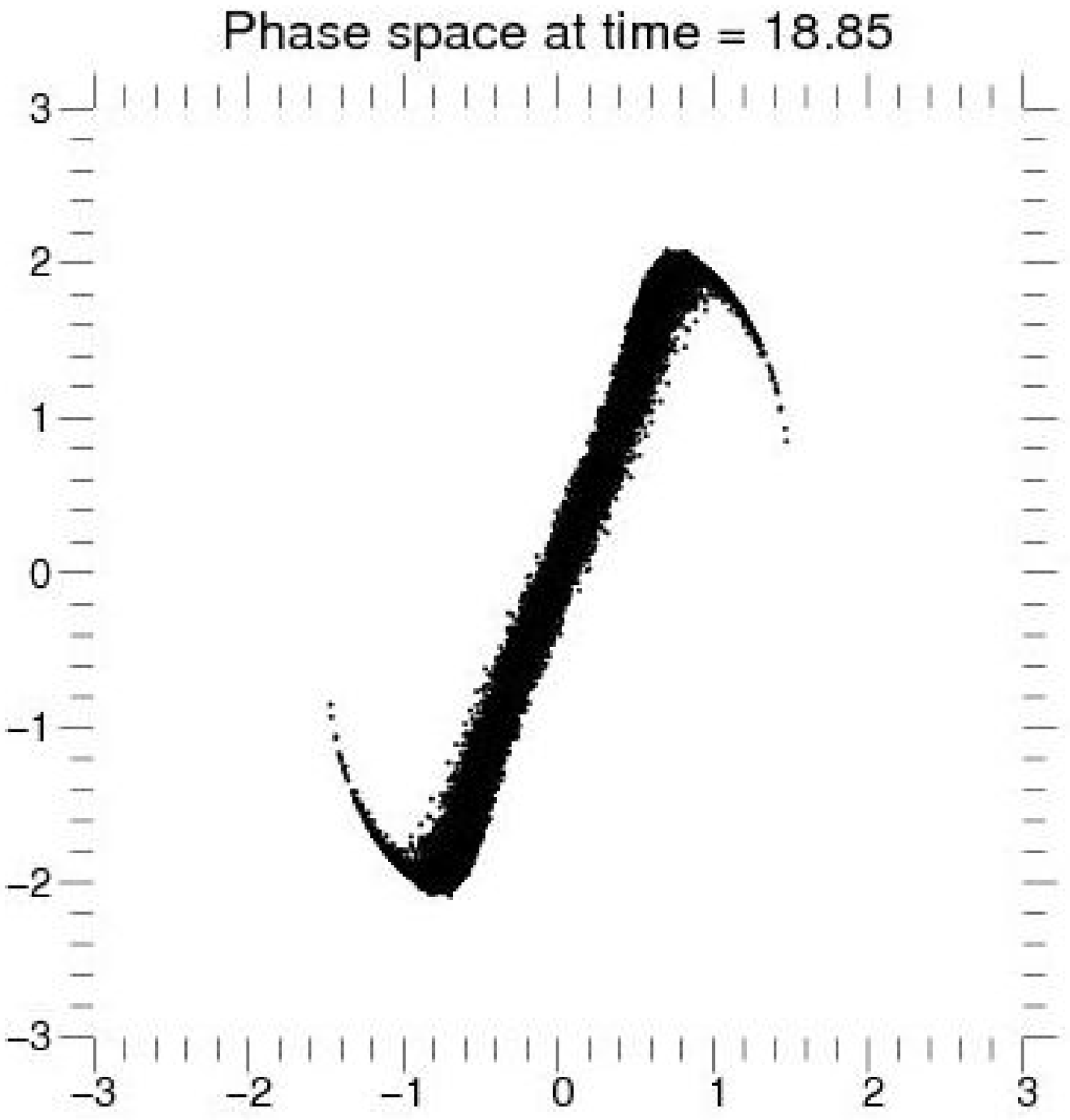}
&
\includegraphics[height=4cm]{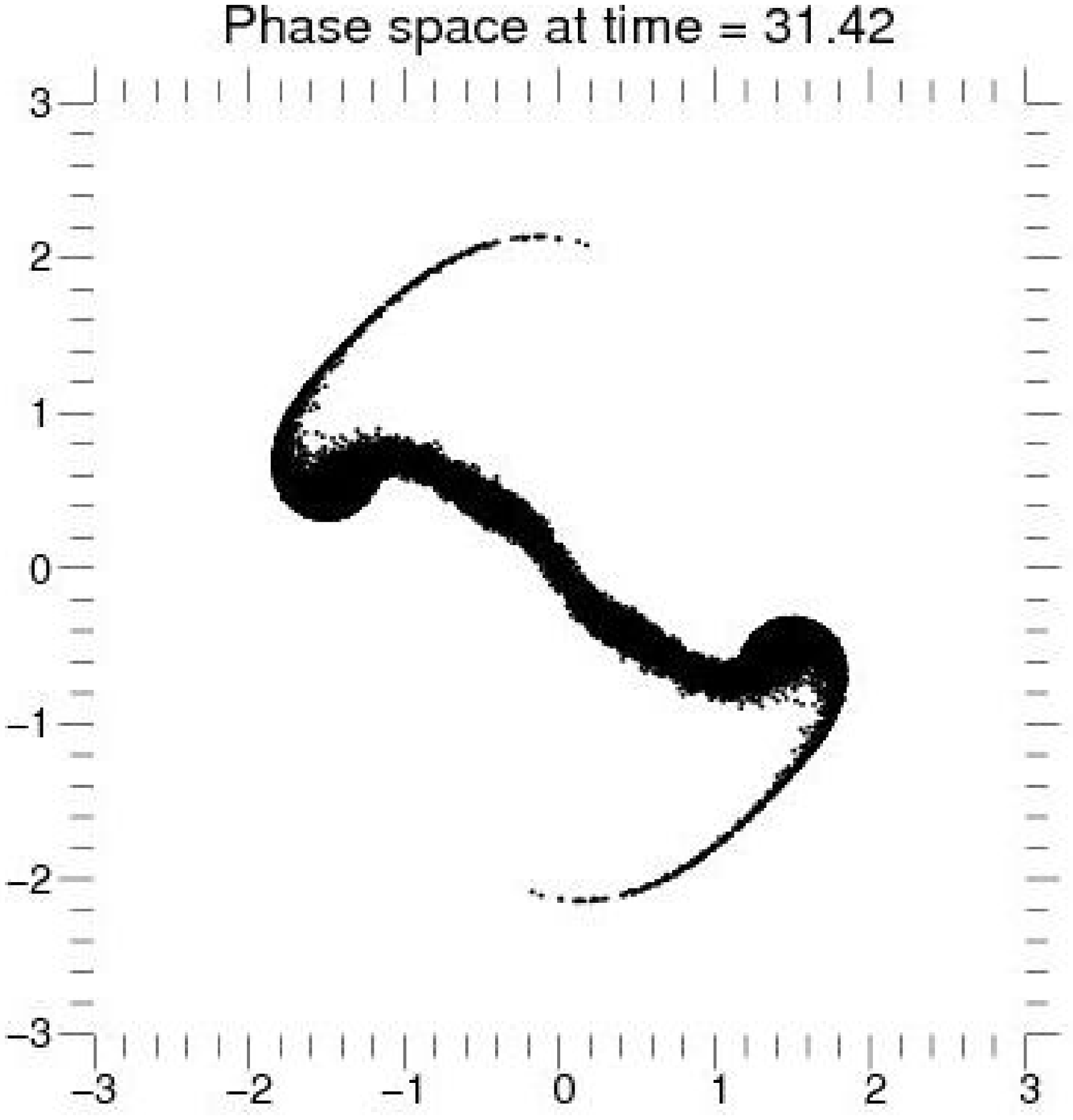}
&
\includegraphics[height=4cm]{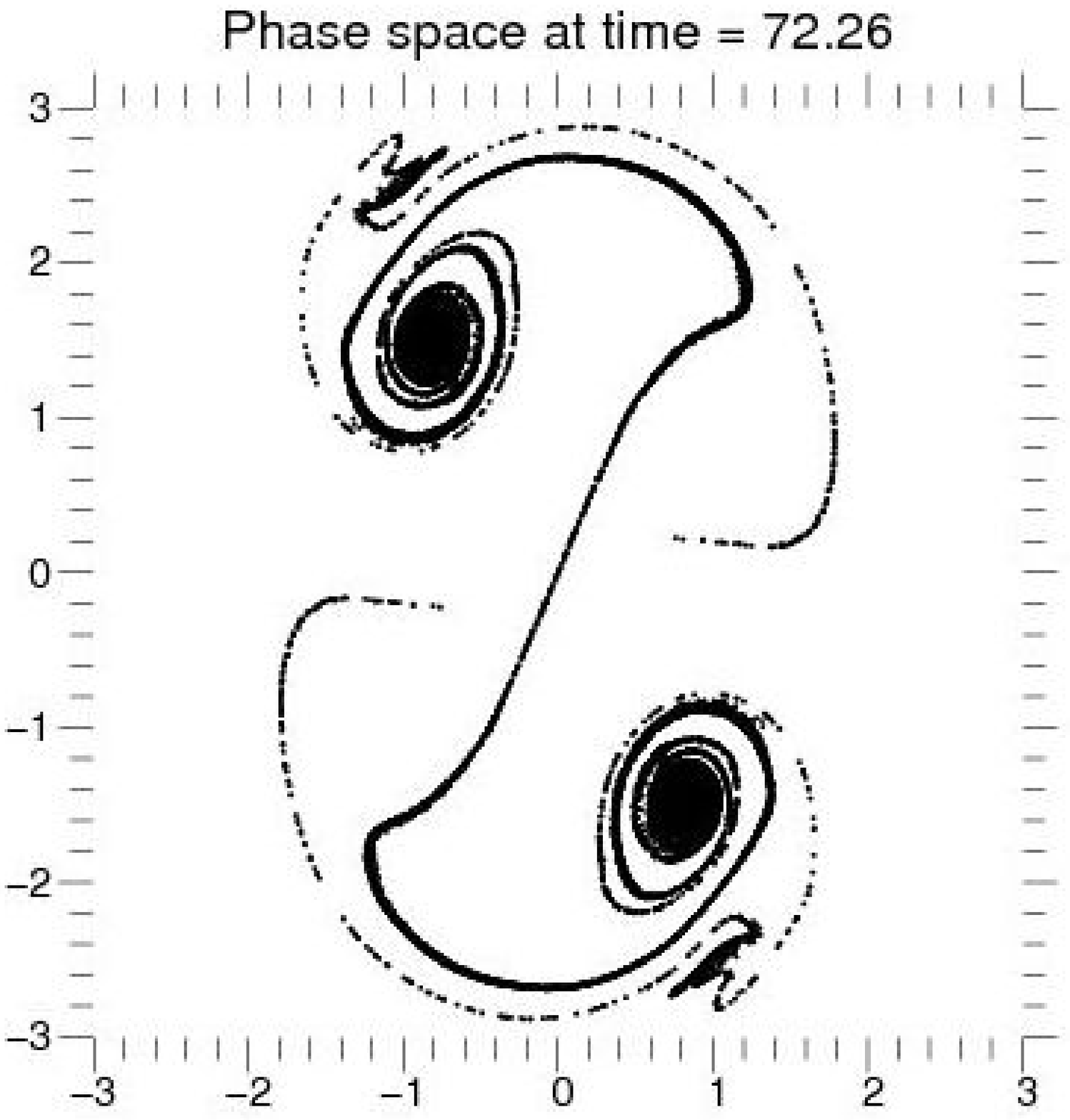}
\\ \\
\includegraphics[height=4cm]{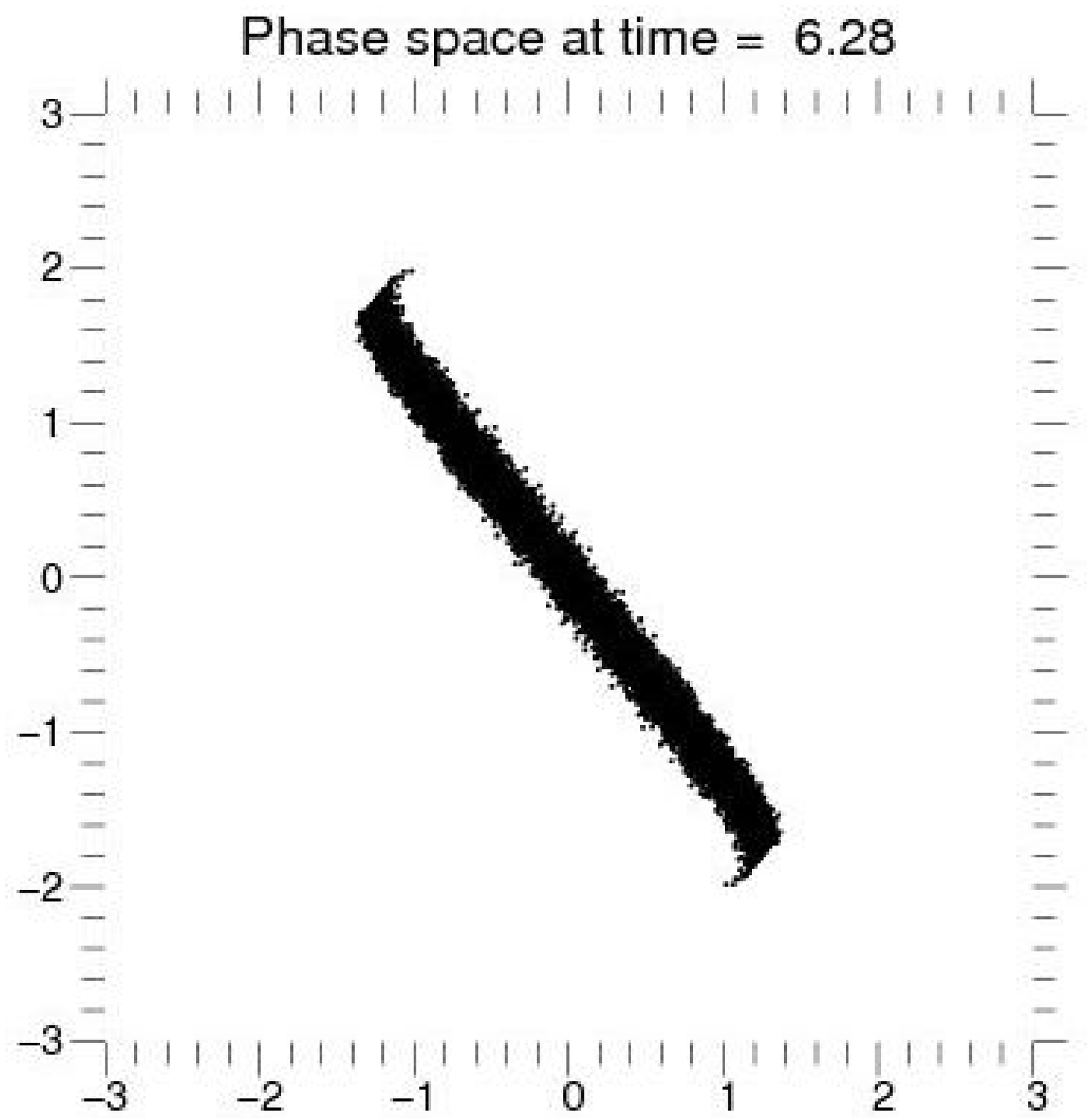}
&
\includegraphics[height=4cm]{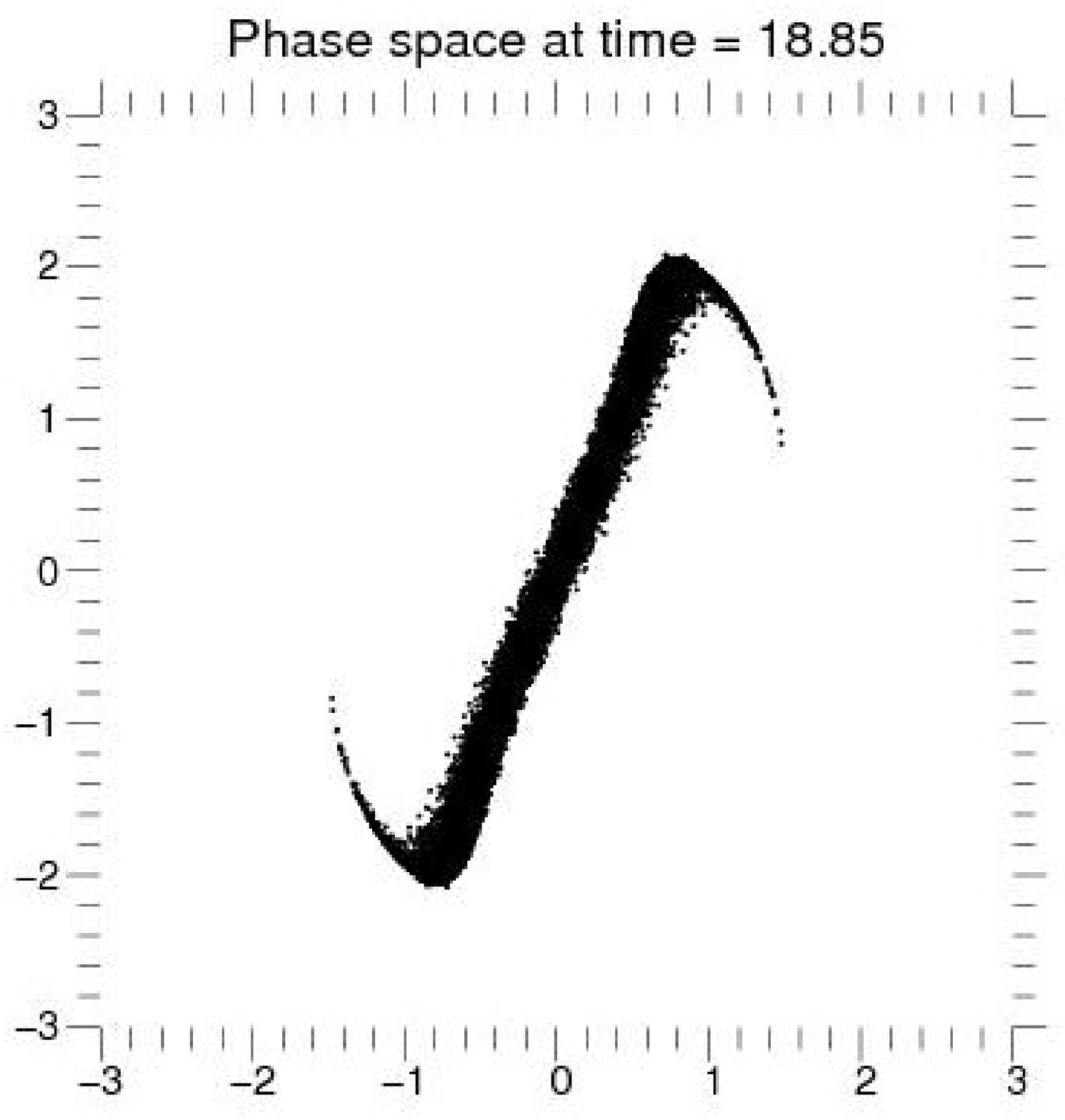}
&
\includegraphics[height=4cm]{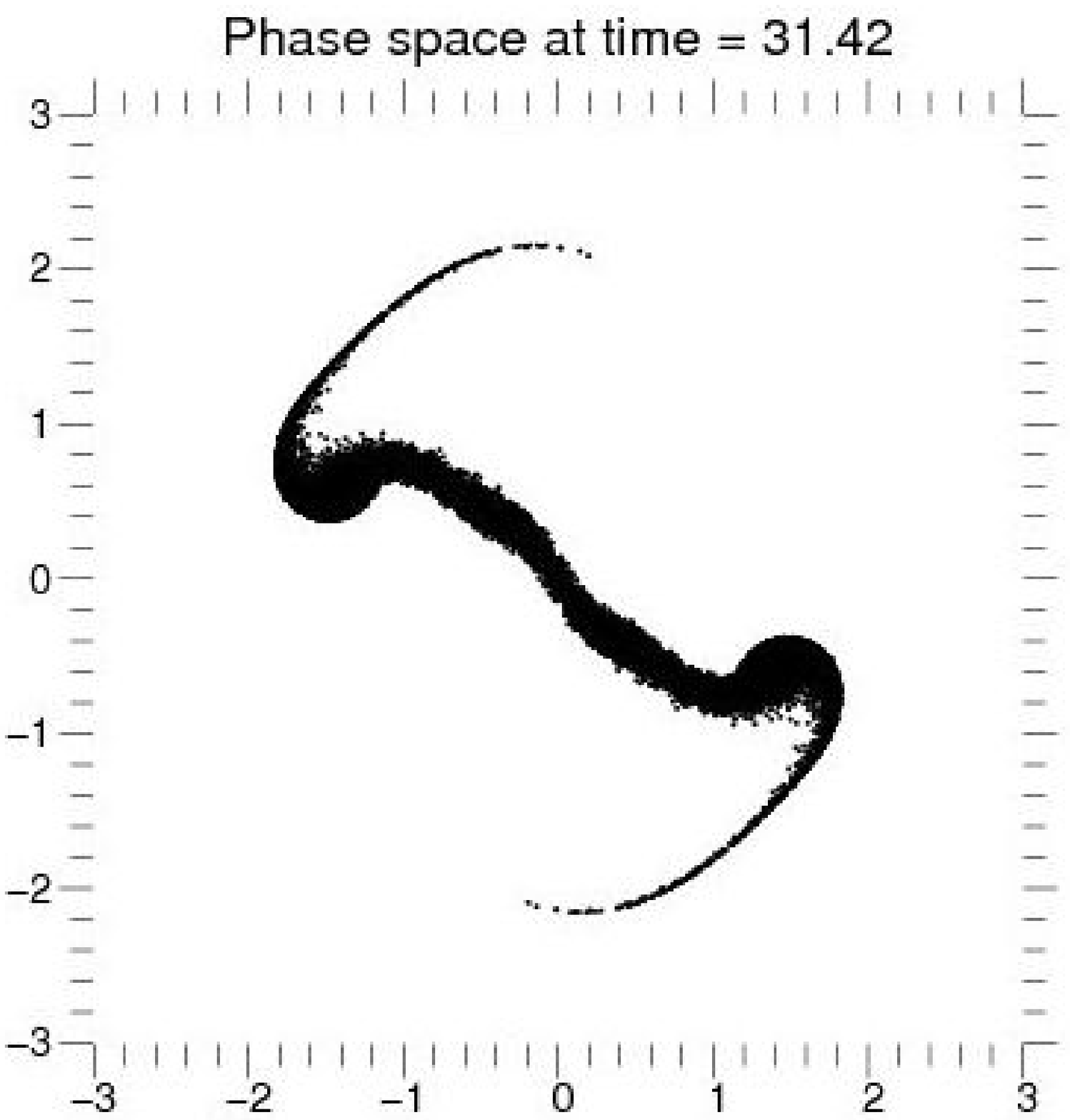}
&
\includegraphics[height=4cm]{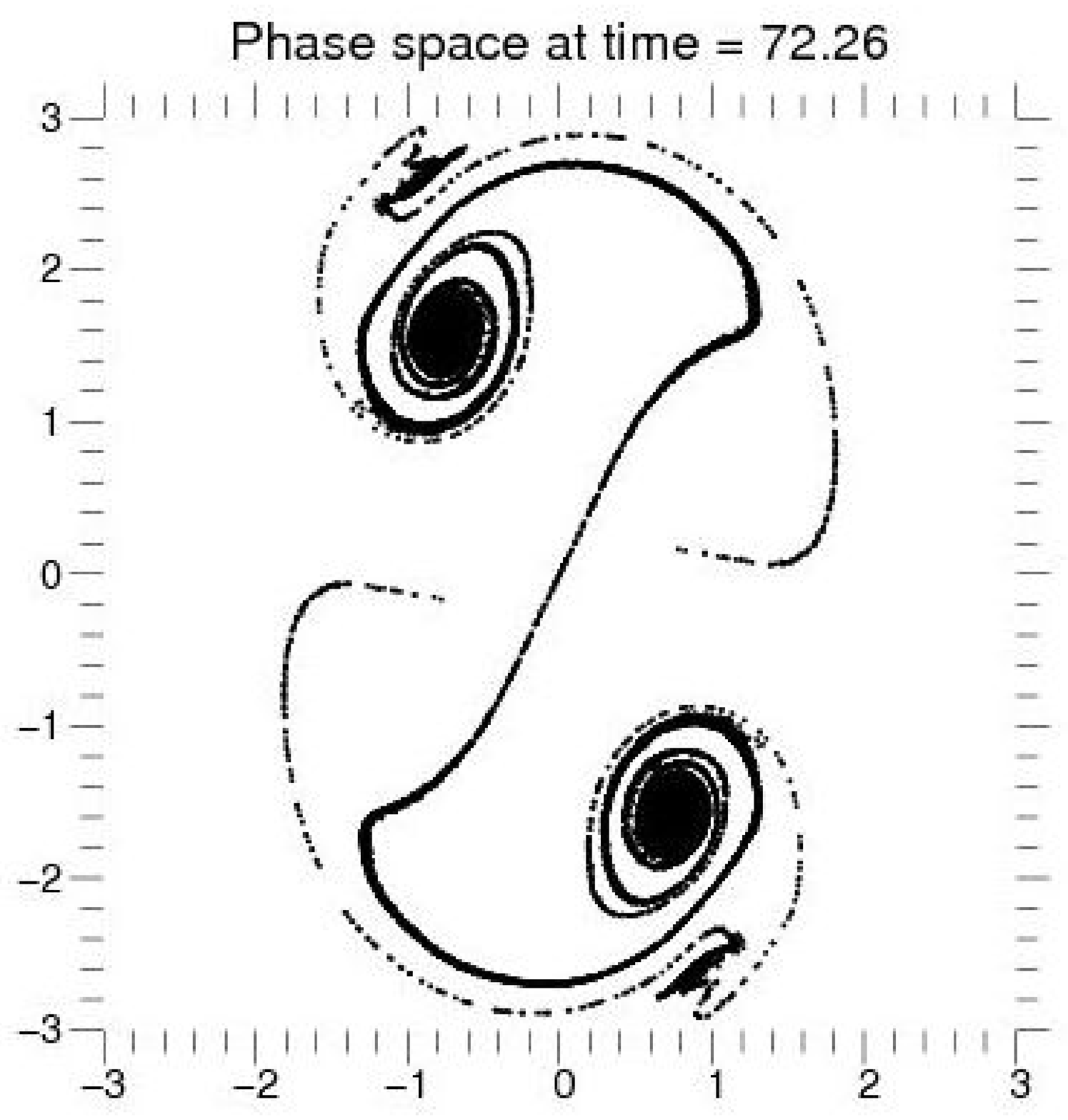}
\end{tabular}
\caption{Left: beam at time 6.26, center left: 
beam at time 18.85, center right : beam at time 31.42, right : beam at time 72.26. 
Top : Simulation provided with the usual PIC method,
bottom:  Simulation provided with the two-scale PIC method.}\label{fig5}
\end{figure}

\medskip
The case when 
\begin{equation}
H_1(\omega_1 \tau) = \cos^2(\tau),
\end{equation}
generates a real effect on the long term.
In other words, the terms containing $H_1$ in (\ref{Zut5bis}) are
\begin{gather}
\frac{1}{2\pi} \int_0^{2\pi}\sin(\sigma) \cos^2(\sigma)\big(\cos(\sigma)q+\sin(\sigma)u_r
\big)\, d\sigma =  -\frac{1}{8}u_r,
\\
\frac{1}{2\pi} \int_0^{2\pi} \cos(\sigma)\cos^2(\sigma)\big(\cos(\sigma)q+\sin(\sigma)u_r
\big)\, d\sigma =  \frac{3}{8}q.
\end{gather}
The simulation results are 
given in figure \ref{fig5}.
The top line shows the time evolution
of the beam simulated with a usual standard PIC method and the bottom line
shows the same simulation with the two-scale PIC method just built.
The simulations coincide with a good degree of accuracy also in this case. 

Finally we consider the case when
\begin{equation}
H_1(\omega_1 \tau) = \cos(2 \tau),
\end{equation}
which has a defocusing effect on the beam as
\begin{gather}
\frac{1}{2\pi} \int_0^{2\pi}\sin(\sigma) \cos(2\sigma)\big(\cos(\sigma)q+\sin(\sigma)u_r
\big)\, d\sigma =  \frac{u_r}{4},
\\
\frac{1}{2\pi} \int_0^{2\pi} \cos(\sigma)\cos(2\sigma)\big(\cos(\sigma)q+\sin(\sigma)u_r
\big)\, d\sigma =  \frac{q}{4}.
\end{gather}

\begin{figure}
\begin{center}
\begin{tabular}{cc}
\includegraphics[height=4.7cm]{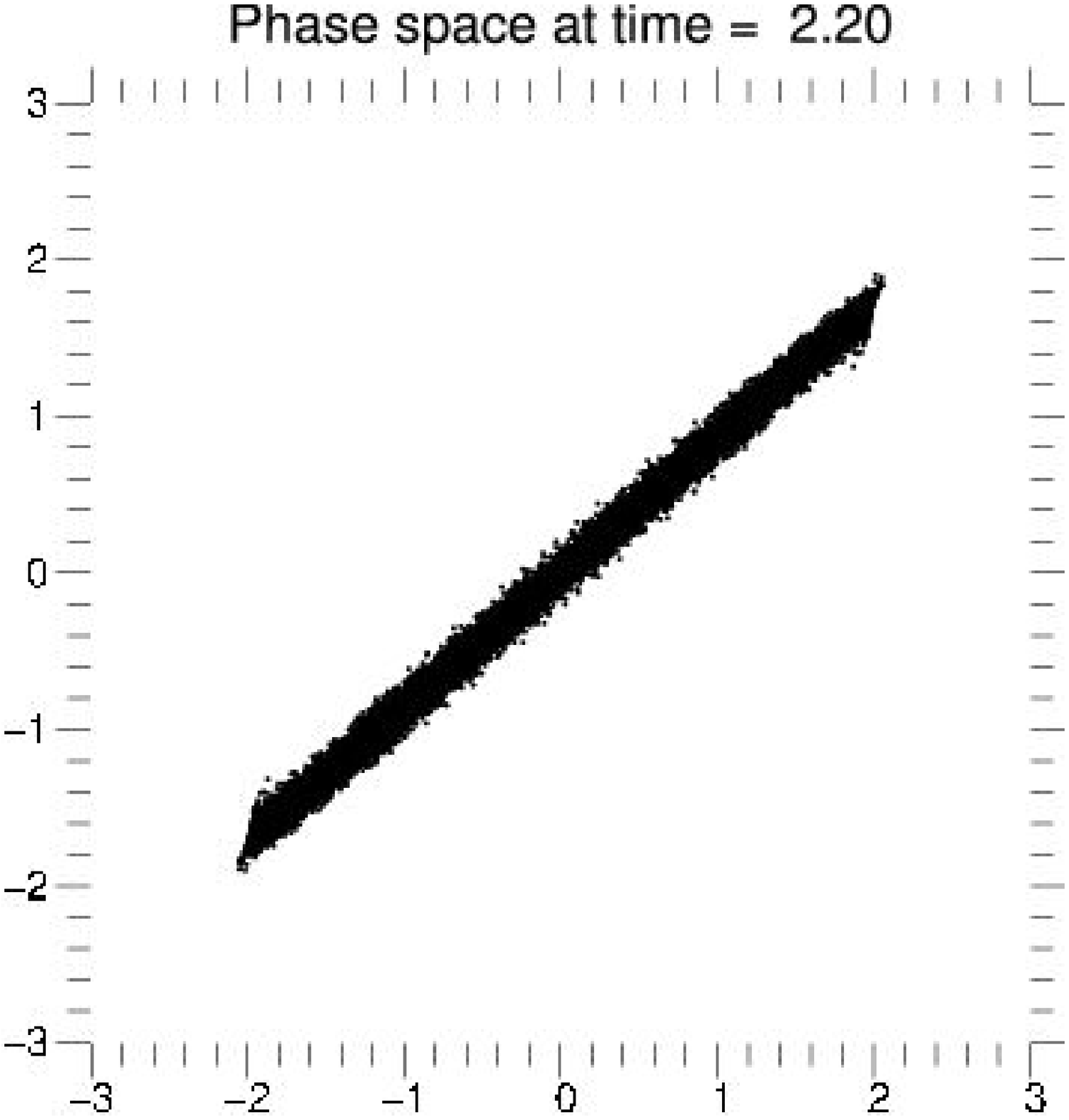}
&
\includegraphics[height=4.7cm]{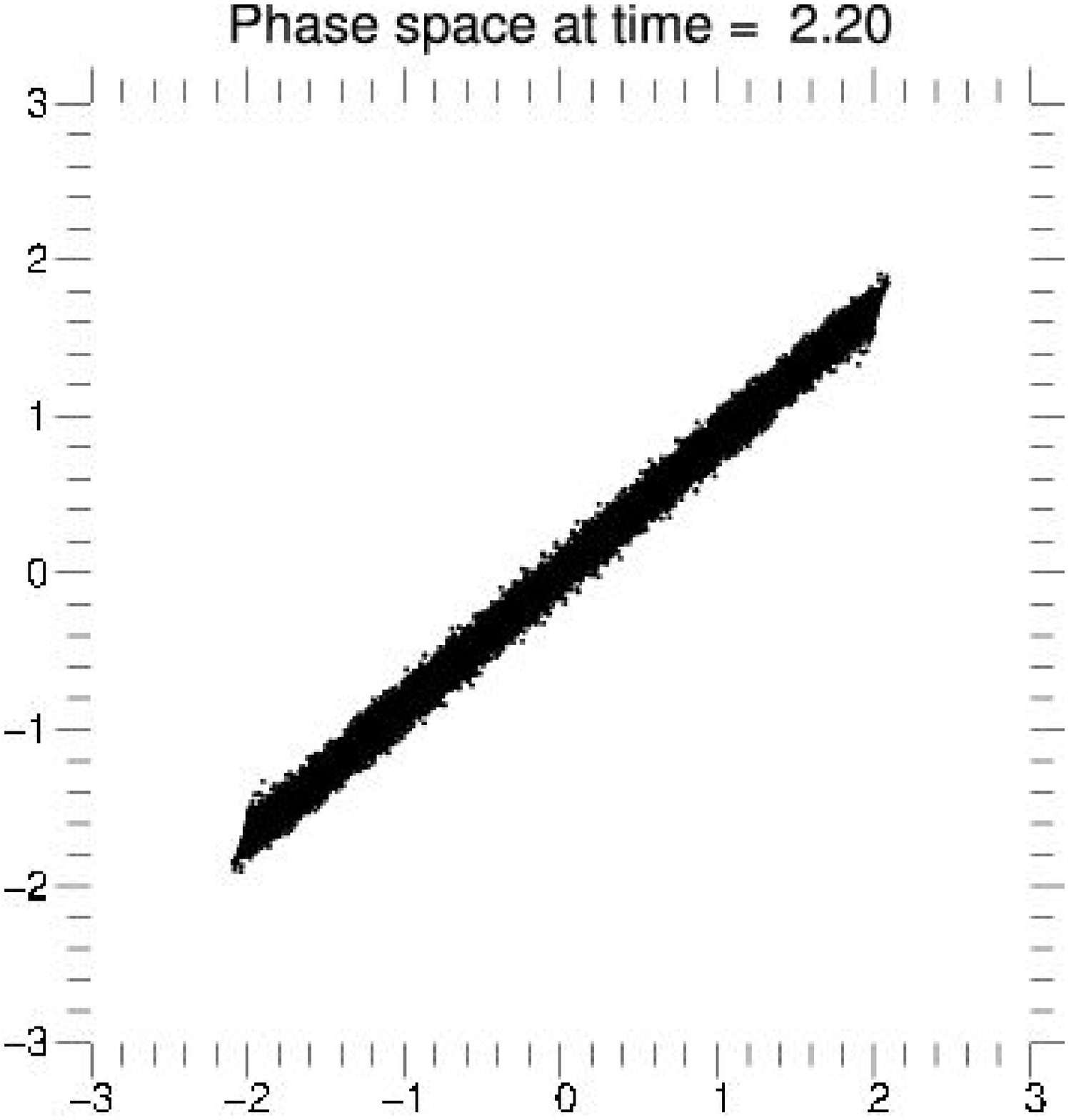}
\end{tabular}
\end{center}
\caption{Beam simulation with an external de-focusing force at time 2.20.
Left: with the usual PIC method. Right: with the two-scale PIC method. }
\label{fig10}
\end{figure}

The result shown in figure \ref{fig10}
shows that also in this case the two-scale PIC method gives
 results that are qualitatively and quantitatively very close to those of the usual PIC method.

\section{Conclusions and perspectives}

We have developed and validated a new PIC solver that can deal very efficiently with a problem involving two time scales. The two-scale solver is based on homogenized equations that have been obtained analytically in our application. For small values of $\epsilon$ the new solver can follow as accurately as the original solver the complex structure of the particle distribution that is generated in a mismatched particle beam
in an accelerator at a small fraction of the cost of the usual PIC solver. These numerical results obtained in 1D  are very promising and the extension to more dimensions as well as the used of semi-Lagrangian or other type of solvers for the Vlasov-like equation instead of a PIC solver deserve to be investigated for accelerator physics as well as for other applications, e.g. strongly magnetized plasmas where the fast motion along the magnetic field lines provides a fast time scale. 
On the other hand, a similar two time scale method can probably also be applied in cases where the homogenized equation can only be computed numerically. 


\vskip0.5cm {\bf Acknowledgments:} The authors acknowledge financial
support from the project HYKE, ``Hyperbolic and
Kinetic Equations: Asymptotics, Numerics, Analysis'' financed by the
European Union under Contract Number 
HPRN-CT-2002-00282. The second author thanks the Universit\' e
Louis Pasteur, where this paper has been partially
written, for hospitality.


\end{document}